\documentclass[12pt]{article}

\usepackage{amsmath,amsthm,amsfonts,amssymb,amscd}
\def\qed{{\unskip\nobreak\hfil\penalty50
\hskip2em\hbox{}\nobreak\hfil$\square$
\parfillskip=0pt \finalhyphendemerits=0\par}\medskip}
\def\prf{\trivlist \item[\hskip \labelsep{\bf Proof\ }]}

\def\l{\langle}
\def\ro{\rho}
\def\ad{{v_0}}
\def\ra{\rangle}

\def\End{{\mathrm {End}}}
\def\Exp{{\mathrm{Exp}}}
\def\Hom{{\mathrm {Hom}}}

\def\col{{\mathrm{col}}}
\def\a{\alpha}

\def\e{\varepsilon}

\def\la{\lambda}

\def\om{\omega}

\def\End{{\mathrm {End}}}
\def\Hom{{\mathrm {Hom}}}

\def\a{\alpha}

\def\e{\varepsilon}

\def\la{\lambda}

\def\ro{{\rho}}

\newtheorem{theorem}{Theorem}[section]
\newtheorem{lemma}[theorem]{Lemma}
\newtheorem{conjecture}[theorem]{Conjecture}
\newtheorem{corollary}[theorem]{Corollary}
\newtheorem{definition}[theorem]{Definition}

\newtheorem{proposition}[theorem]{Proposition}
\newtheorem{remark}[theorem]{Remark}
\newtheorem{example}[theorem]{Example}

\def\Hom{{\mathrm{Hom}}}

\def\res{\!\restriction\!}
\def\A{{\cal A}}

\def\Z{{\mathbb Z}}

\def\Co{{\mathbb C}}

\renewcommand{\qed}{\ \hfill $\blacksquare$}

\newcommand{\bdef}{\begin{definition}}
\newcommand{\blem}{\begin{lemma}}
\newcommand{\bprop}{\begin{proposition}}
\newcommand{\bthm}{\begin{theorem}}
\newcommand{\bcor}{\begin{corollary}}
\newcommand{\bconj}{\begin{conjecture}}
\newcommand{\ede}{\end{definition}}
\newcommand{\elem}{\end{lemma}}
\newcommand{\eprop}{\end{proposition}}
\newcommand{\ethm}{\end{theorem}}
\newcommand{\ecor}{\end{corollary}}
\newcommand{\econj}{\end{conjecture}}
\newcommand{\brem}{\begin{remark}}
\newcommand{\erem}{\end{remark}}

\newcommand{\ba}{\begin{array}}
\newcommand{\ea}{\end{array}}
\newcommand{\bea}{\begin{eqnarray}}

\renewcommand{\mod}{\mbox{mod}}

\title{\huge On representing some  lattices as lattices of intermediate subfactors of finite index\\}
\author{
{\sc Feng Xu}\footnote{Supported in part by NSF.}\\
Department of Mathematics\\
University of California at Riverside\\
Riverside, CA 92521\\
E-mail: {\tt xufeng@math.ucr.edu}}
\begin{document}
\date{}
\maketitle

\begin{abstract}
We prove that the very simple lattices which consist of a largest, a
smallest and $2n$ pairwise incomparable elements where $n$ is a
positive integer can be realized as the lattices of intermediate
subfactors of finite index and finite depth. Using the same
techniques, we  give a necessary and sufficient condition for
subfactors coming from Loop groups of type $A$ at generic levels to
be maximal.

\end{abstract}

\newpage

\section{Introduction}
Let $M$ be a factor and $N$ a subfactor of $M$ which is irreducible,
i.e.,$N'\cap M= \Co$. Let $K$ be an intermediate von Neumann
subalgebra for the inclusion $N\subset M.$ Note that $K'\cap
K\subset N'\cap M = \Co,$ $K$ is automatically a factor. Hence the
set of all intermediate subfactors for $N\subset M$ forms a lattice
under two natural operations $\wedge$ and $\vee$ defined by:
\[
K_1\wedge K_2= K_1\cap K_2, K_1\vee K_2= (K_1\cup K_2)''.
\]
Let $G_1$ be a group and $G_2$ be a subgroup of $G_1$.  An interval
sublattice $[G_1/G_2]$ is the lattice formed by all intermediate
subgroups $K, G_2\subseteq K\subseteq G_1.$

By cross product construction and Galois correspondence,  every
interval sublattice of finite groups can be realized as intermediate
subfactor lattice of finite index. The study of intermediate
subfactors has been very active in recent years (cf.
\cite{BJ},\cite{GJ}, \cite{ILP}, \cite{Longo4},  \cite{Wat} and
\cite{ SW} for only a partial list). By a result of S. Popa (cf.
\cite{Pol}), if a subfactor $N\subset M$ is irreducible and has
finite index, then the set of intermediate subfactors between $N$
and  $M$ is finite. This result was also independently proved by Y.
Watatani (cf. \cite{Wat}). In \cite{Wat}, Y. Watatani investigated
the question of which finite lattices can be realized as
intermediate subfactor lattices. Related questions were further
studied by P. Grossman and V. F. R. Jones in \cite{GJ} under certain
conditions. As emphasized in \cite{GJ}, even for a lattice which
shapes like a Hexagon and consists of six elements, it is not clear
if it can be realized as intermediate subfactor lattice with finite
index. This question has been solved recently by M. Aschbach in
\cite{Asch} among other things. In \cite{Asch}, M. Aschbach
constructed a finite group $G_1$ with a subgroup $G_2$ such that the
interval sublattice $[G_1/G_2]$ is a Hexagon. The lattices that
appear in \cite{GJ,Wat,Asch} can all be realized as interval
sublattice of finite groups.
\par It turns out that which
finite lattice can be realized as an interval sublattice $[G_1/G_2]$
with $G_1$ finite is an old problem in finite group theory. See
\cite{Palfy} for an excellent review and a list of references.

Most of the attention has been focused on the very simple lattice
$M_n$ consisting of a largest, a smallest and $n$ pairwise
incomparable elements. For $n=1, 2, q+1$ (where $q$ is a prime
power), examples of $M_n$ have been found in the finite solvable
groups. After the first interesting examples found by W. Feit (cf.
\cite{Feit}), A. Lucchini (cf. \cite{Lucc})  discovered new series
of examples for $n=q+2$ and for $n=\frac{(q^t+1)}{(q+1)}+1$ where
$t$ is an odd prime.\par At the present the only values of $n$ for
which $M_n$ occurs as an interval sublattice of a finite group are
$n=1,2, q+1,q+2, \frac{(q^t+1)}{(q+1)}+1 $  where $t$ is an odd
prime. The smallest undecided case is $n=16.$ In a major
contribution to the problem about subgroup lattices of finite groups
in \cite{Balu}, R. Baddeley and A. Lucchini have reduced the problem
of realizing $M_n$ as interval sublattice of finite groups to a
collection of questions about finite simple groups. These questions
are still quite hard, but eventually they might be resolved using
the classification of finite simple groups. In this paper, the
authors are cautious, but their ultimate goal seems to be to show
that the list above is complete. In view of the above results about
finite groups, it seems an interesting problem to ask if $M_{16}$
can be realized as the lattice of intermediate subfactors with
finite index.  This problem is the main motivation for our paper.
One of the main results of this paper, Theorem \ref{mn}, states that
all $M_{2n}$ are realized as the lattice of intermediate subfactors
of a pair of hyperfinite type $III_1$ factors with  finite depth.
Note that by \cite{SW} this implies that $M_{2n}$ can also be
realized as the lattice of intermediate subfactors of a pair of
hyperfinite type $II_1$ factors with  finite depth. Thus modulo the
conjectures of R. Baddeley, A. Lucchini and possibly others we have
an infinite series of lattices which can be realized by the lattice
of intermediate subfactors with finite index and finite depth but
can not be realized by interval sublattices of finite groups.\par
The subfactors which realize $M_{2n}$ are ``orbifold subfactors" of
\cite{EK, BE2, Xo}, and we are motivated to examine these subfactors
by the example of lattice of type $M_6$ in \cite{GJ} which is
closely related to an ${\Bbb Z}_2$ orbifold. To explain their
construction, after first two preliminary sections, we will first
review the  result of A. Wassermann (cf. \cite{Jh1},\cite{Wass})
about Jones-Wassermann subfactors (cf. Remark (\ref{jxx})) coming
from representations of Loop groups of type $A$ in section
\ref{typea}. Section \ref{orb} is then devoted to a description of
``orbifold subfactors" from an induction point of view. Although it
is not too hard to show that the subfactor contains $2n$
incomparable intermediate subfactors, the hard part of the proof of
Theorem \ref{mn} is to show that there are no more intermediate
subfactors. Here we give a brief explanation of basic ideas behind
our proof and describe how the paper is structured.  We will use
freely notations and concepts that can be found in preliminary
sections. Let $\rho(M)\subset M$ be a subfactor and $M_1$ be an
intermediate subfactor. In our examples below all factors are
isomorphic to the hyperfinite type $III_1$ factor, and
$\rho\bar\rho$ are direct sums of sectors from a  set $\Delta$ with
finitely many irreducible sectors and a non-degenerate braiding.
Here we use the endomorphism theory pioneered by R. Longo (cf.
\cite{Longo1}). Since $M_1$ is isomorphic to $M,$ we can choose an
isomorphism $c_1: M\rightarrow M_1.$ Denote by $c_2= c_1^{-1}\rho$
we have $\rho = c_1c_2$ where $c_1, c_2\in {\End(M)}.$ Note that
$c_1\bar{c_1}\prec \rho\bar\rho$ is in $\Delta.$ Our basic idea to
investigate the property of $c_1$ is to consider the following set
$H_{c_1}:= \{ [a] | a\prec \lambda c_1, \lambda\in \Delta, \ a \ \
{\mathrm irreducible} \}.$ Since $\Delta$ has finitely many
irreducible sectors, $H_{c_1}$ is a finite set. Moreover since
$c_1\bar c_1\in \Delta,$ an induction method using braidings as in
\cite{Xb} is available. This induction method was used by the author
in \cite{Xb} to study subfactors from conformal inclusions, and
developed further by J. B\"{o}ckenhauer-D. Evans and J.
B\"{o}ckenhauer-D. Evans-Y. Kawahigashi in \cite{BE1}, \cite{BE2},
\cite{BE3} , \cite{BE4}, \cite{BEK1}, and \cite{BEK2}, and leads to
strong constraints on the set $H_{c_1}.$ Thus by using $\lambda \in
\Delta$ to act from the left on $c_1,$ one may hope to find what
$c_1$ is made of.  In the cases of Theorem \ref{mn} and Cor.
\ref{max}, it turns out that there is a sector $c$ in $H_{c_1}$ with
smallest index such that $c_1\prec \lambda c,$ and $c$ is close to
be an automorphism ( It is an automorphism in the case of Cor.
\ref{max}), and the corresponding subfactors have been well studied
as those in \cite{Xb}. In the simplest case $n=2,$ due to $A-D-E$
classification of graphs with norm less than $2,$ the above idea can
be applied directly to give a rather quick proof of Th. \ref{mn}. We
refer the reader to the paragraph after Th. \ref{mn} which
illustrates the above idea. \par When $n>2,$ the norms of fusion
graphs are in general greater than $2,$ no $A-D-E$ classification is
available, and this is the main problem we must resolve to carry out
the above idea. To explain our method, we note that $S$ matrix as
defined in equation (\ref{Smatrix}) has the property
\[
|\frac{S_{\lambda\mu}}{S_{1\mu}}|\leq \frac{S_{\lambda 1}}{S_{1 1}}
, \forall \mu\] and
\[
|\frac{S_{\lambda\mu}}{S_{1\mu}}|= \frac{S_{\lambda 1}}{S_{1 1}} ,
\forall \mu\] iff $\lambda$ is an automorphism, i.e., $\lambda$ has
the smallest index $1.$  Our first observation is that for small
index (close to $1$) sectors $c,$ certain entries of $S$-matrix like
quantities (cf. Def. \ref{psiw}, Cor. \ref{smallc}) called $\psi$-
matrix attain their extremum just like $S$-matrices. Hence to detect
these small index sectors, we need to have a good estimation of the
entries of $\psi$- matrix. In view of the Verlinde formula (cf.
equation \ref{Verlinde}) relating $S$-matrix with fusion rules, it
is natural to use the known fusion rules to estimate $\psi$ matrix.
However, since the definition of $\psi$ involves sectors which are
not braided, the above idea does not work  unless one can show that
certain intertwining operators are central (cf. Th. \ref{centu} and
Section \ref{centrality} for discussions) . Our second observation
is that a class of intertwining operators in Definition \ref{uwm}
are central (cf. Th. \ref{centu}). Thanks to a number of known
results about representations of Loop groups of type $A,$  we show
that the assumption of Th. \ref{centu} is verified in our case(cf.
Prop. \ref{lemh}).

This  allows us to show that certain sector with small index does
exist (cf. Cor. \ref{smallc}), we can indeed find that $c_1$ is made
of known subfactors. After a straightforward calculations involving
known fusion rules in Prop. \ref{nonprime}, we are able to finish
the proof of Th. \ref{mn} for general $n$. \par In the last section
we discuss a few related issues. Conjecture \ref{conj1} is
formulated which is equivalent to centrality of certain intertwining
operators (cf. Prop. \ref{c1=c}), and this is motivated by our proof
of Th. \ref{mn}. We show in Prop. \ref{ub} that these intertwining
operators are central on a subspace which is a linear span of
products of (cf. Def. \ref{cap}) cups, caps and braiding operators
only. These motivate us to make Conjecture \ref{conj2} which claims
that the subspace is fact the whole space. In view of recent
development using category theory (cf. \cite{FSS}), both conjectures
can in fact be stated in categorical terms, and we don't know any
counter examples in the categorical setting. In Prop. \ref{c3=c1} we
prove that a weaker version of Conjecture \ref{conj2} implies
Conjecture \ref{conj1}, and from this we are able to prove
Conjecture \ref{conj1} for modular tensor category from $SU(n)$ at
level $k$ (cf. Cor. \ref{hecke}).\par In section \ref{maxsection} we
give applications of Cor. \ref{hecke}. To explain these
applications, recall that a subfactor $N\subset M$ is called maximal
if $M_1$ is an intermediate von Neumann algebra between $N$ and $M$
implies that $M_1=M$ or $M_1=N.$ This notion is due to V. F. R.
Jones when he outlined an interesting programme to investigate
questions in group theory using subfactors (cf. \cite{Jh}). In the
case when $M$ is the crossed product of $N$ by a finite group $G,$
it is easy to see that $N\subset M$ is maximal iff $G$ is an abelian
group of prime order. Hence for most of $G$ the corresponding
subfactor is not maximal. Cor. \ref{max} gives a necessary and
sufficient condition for subfactors from representations $\lambda$
of $SU(n)$ at level $k\neq n\pm2, n$ to be maximal: $\lambda$ is
maximal iff $\lambda$ is not fixed by a nontrivial cyclic
automorphism of extended Dynkin diagram (Such cyclic automorphisms
generate a group isomorphic to ${\Bbb Z}_n$). Hence it follows from
Cor. \ref{max} that most of such $\lambda$ are maximal. For an
example, if $k\neq n\pm2, n$, $k$ and $n$ are relatively prime, then
all $\lambda$ are maximal. \par

Besides Propositions and Theorems that have been already mentioned,
The first two preliminary sections are about sectors, covariant
representations, braiding-fusion equations, Yang-Baxter equations,
Rehren's $S,T$ matrices. The third preliminary section summarizes
properties of an induction method from \cite{Xb}. These properties
have been extensively studied and applied in subsequent work in
\cite{BE1},\cite{BE2},\cite{BE3},\cite{BE4}, \cite{BEK1} and
\cite{BEK2} from a different point of view where induction takes
place between two different but isomorphic algebras, and we recall a
dictionary relating these two as provided in \cite{Xj}. We think
that  in this paper it is simpler to take the point of view of
\cite{Xb} when discussing intermediate subfactors, and it is
convenient to represent these intermediate subfactors as the range
of endomorphisms of one fixed factor,  so we do not have to switch
between different but isomorphic algebras.

Using the dictionary we translate some properties of relative
braidings and local extensions from \cite{BE4} to our setting (cf.
Prop. \ref{loclr}). The next two preliminary sections are devoted to
subfactors from representations of $SU(n)$ at level $k$ and its
extensions. We collect a few properties about fusion rules, $S$
matrices, and we define the subfactor which appears in Th. \ref{mn}.
In Prop. \ref{2n} we show that this subfactor contains $2n$
incomparable proper intermediate subfactors. \par

The author would like to thank Professor M. Aschbacher for useful
discussions on subgroup lattices of finite groups, and especially
Professor V. F. R. Jones for helpful suggestions and encouragement.

\section{Preliminaries}

For the convenience of the reader we collect here some basic notions
that appear in this paper. This is only a guideline and the reader
should look at the references such as preliminary sections of
\cite{KLX} for a more complete treatment.

\subsection{Sectors}
Let $M$ be a properly infinite factor and  $\text{\rm End}(M)$ the
semigroup of
 unit preserving endomorphisms of $M$.  In this paper $M$ will always
be the unique hyperfinite $III_1$ factors. Let $\text{\rm Sect}(M)$
denote the quotient of $\text{\rm End}(M)$ modulo unitary
equivalence in $M$. We  denote by $[\rho]$ the image of $\rho \in
\text{\rm End}(M)$ in  $\text{\rm Sect}(M)$.\par
 It follows from
\cite{Longo1} and \cite{Longo2} that $\text{\rm Sect}(M)$, with $M$
a properly infinite  von Neumann algebra, is endowed with a natural
involution $\theta \rightarrow \bar \theta $  ; moreover, $\text{\rm
Sect}(M)$ is
 a semiring.\par
Let $\rho \in \text{\rm End}(M)$ be a normal faithful conditional
expectation $\epsilon: M\rightarrow \rho(M)$.  We define a number
$d_\epsilon$ (possibly $\infty$) by:
$$
d_\epsilon^{-2} :=\text{\rm Max} \{ \lambda \in [0, +\infty)|
\epsilon (m_+) \geq \lambda m_+, \forall m_+ \in M_+ \}$$ (cf.
[PP]).\par
 We define
$$
d = \text{\rm Min}_\epsilon \{ d_\epsilon |  d_\epsilon < \infty \}.
$$   $d$ is called the statistical dimension of  $\rho$ and $d^2$ is called the Jones index of $\rho$.
It is clear from the definition that  the statistical dimension  of
$\rho$ depends only on the unitary equivalence classes  of  $\rho$.
The properties of the statistical dimension can be found in
\cite{Longo1},\cite{Longo2} and \cite{Longo3}.\par Denote by
$\text{\rm Sect}_0(M)$ those elements of $\text{\rm Sect}(M)$ with
finite statistical dimensions. For $\lambda $, $\mu \in \text{\rm
Sect}_0(M)$, let $\text{\rm Hom}(\lambda , \mu )$ denote the space
of intertwiners from $\lambda $ to $\mu $, i.e. $a\in \text{\rm
Hom}(\lambda , \mu )$ iff $a \lambda (x) = \mu (x) a $ for any $x
\in M$. $\text{\rm Hom}(\lambda , \mu )$  is a finite dimensional
vector space and we use $\langle \lambda , \mu \rangle$ to denote
the dimension of this space. $\langle  \lambda , \mu \rangle$
depends only on $[\lambda ]$ and $[\mu ]$. Moreover we have $\langle
\nu \lambda , \mu \rangle = \langle \lambda , \bar \nu \mu \rangle
$, $\langle \nu \lambda , \mu \rangle = \langle \nu , \mu \bar
\lambda \rangle $ which follows from Frobenius duality (See
\cite{Longo2} ).  We will also use the following notation: if $\mu $
is a subsector of $\lambda $, we will write as $\mu \prec \lambda $
or $\lambda \succ \mu $.  A sector is said to be irreducible if it
has only one subsector. \par For any $\rho\in \text{\rm End}(M)$
with finite index, there is a unique standard minimal inverse
$\phi_\rho: M\rightarrow M$ which satisfies
\[
\phi_\rho(\rho(m)m'\rho(m''))= m\phi_\rho(m')m'', m,m',m''\in M.
\]
$\phi_\ro$ is completely positive.  If $t\in \Hom(\rho_1,\rho_2)$
then we have
\begin{equation}\label{phiprop}
d_{\rho_1} \phi_{\rho_1} (mt)= d_{\rho_2} \phi_{\rho_2} (tm), m\in M
\end{equation}
\subsection{Sectors from conformal nets and their representations}
We refer the reader to \S3 of \cite{KLX} for definitions of
conformal nets and their representations. Suppose a conformal net
$\A$ and a representation $\lambda$ is given. Fix an open interval I
of the circle and Let $M:=\A(I)$ be a fixed type $III_1$ factor.
Then $\lambda$ give rises to an endomorphism still denoted by
$\lambda$ of $M$. We will recall some of the results of \cite{R2}
and introduce notations. \par Suppose $\{[\lambda] \}$ is a finite
set  of all equivalence classes of irreducible, covariant,
finite-index representations of an irreducible local conformal net
$\A$. We will use $\Delta_\A$ to denote all finite index
representations of net $\A$ and will use the same notation
$\Delta_\A$ to denote the corresponding sectors of $M$\footnote{Many
statements in this section and section \ref{inductionsection} hold
true in general case when  the the set $\{[\lambda] \}$ is only
braided (cf. \cite{BEK1}) and we hope to consider such cases
elsewhere.}.

We will denote the conjugate of $[\lambda]$ by $[{\bar \lambda}]$
and identity sector (corresponding to the vacuum representation) by
$[1]$ if no confusion arises, and let $N_{\lambda\mu}^\nu = \langle
[\lambda][\mu], [\nu]\rangle $. Here $\langle \mu,\nu\rangle$
denotes the dimension of the space of intertwiners from $\mu$ to
$\nu$ (denoted by $\text {\rm Hom}(\mu,\nu)$).  We will denote by
$\{T_e\}$ a basis of isometries in $\text {\rm
Hom}(\nu,\lambda\mu)$. The univalence of $\lambda$ and the
statistical dimension of (cf. \S2  of \cite{GL1}) will be denoted by
$\omega_{\lambda}$ and $d{(\lambda)}$ (or $d_{\lambda})$)
respectively. The unitary braiding operator $\epsilon (\mu,
\lambda)$ (cf. \cite{GL1} ) verifies the following
\begin{proposition}\label{bfe}
\text{\rm (1)}  Yang-Baxter-Equation (YBE)
$$
\e(\mu , \gamma ) \mu (\e(\lambda , \gamma )) \e(\lambda , \mu ) =
\gamma (\e(\lambda , \mu )) \e(\lambda , \gamma )\lambda (\e(\mu ,
\gamma ))\, .
$$

\text{\rm (2)}  Braiding-Fusion-Equation (BFE)

For any $w\in \text{\rm Hom} (\mu \gamma , \delta )$
\begin{align*} \e(\lambda , \delta )\lambda (w) = w\mu (\e(\lambda , \gamma ))
\e(\lambda , \mu ) \\
\e(\delta , \lambda )w = \lambda (w) \e(\mu , \lambda ) \mu
(\e(\gamma , \lambda )) \, \\
\e( \delta,\lambda )^* \lambda (w) = w\mu (\e( \gamma,\lambda )^*)
\e(\mu,\lambda )^* \\
\e(\lambda, \delta )^* \lambda (w) = w\mu (\e( \gamma,\lambda )^*)
\e(\lambda,\mu )^*
\end{align*}
\end{proposition}

\blem\label{monodromy} If $\lambda,\mu$ are irreducible, and
$t_\nu\in \Hom(\nu, \lambda\mu)$ is an isometry, then $t_\nu
\e(\mu,\lambda)\e(\lambda,\mu) t_\nu^* =
\frac{\om_\nu}{\om_\lambda\om_\mu}.$ \elem

By Prop. \ref{bfe}, it follows that if $t_i\in \Hom (\mu_i,\lambda)$
is an isometry, then
\[
\e(\mu,\mu_i)\e(\mu_i,\mu)=t_i^*\e(\mu,\la)\e(\lambda,\mu) t_i
\]
We shall always identify the center of $M$ with $\Co.$  Then we have
the following

\blem\label{escalar} If
\[\e(\mu,\la)\e(\lambda,\mu)\in \Co,\]
then \[\e(\mu,\mu_i)\e(\mu_i,\mu)\in \Co, \forall \mu_i\prec \la. \]
\elem

Let $\phi_\lambda$ be the unique minimal left inverse of $\lambda$,
define:
\begin{equation}\label{ymatrix}
Y_{\lambda\mu}:= d(\lambda)  d(\mu) \phi_\mu (\epsilon (\mu, \lambda)^*
\epsilon (\lambda, \mu)^*),
\end{equation}
where $\epsilon (\mu, \lambda)$ is the unitary braiding operator
 (cf. \cite{GL1} ). \par
We list two properties of $Y_{\lambda \mu}$ (cf. (5.13), (5.14) of
\cite{R2}):
\begin{lemma}\label{yprop}
\begin{equation*}
Y_{\lambda\mu} = Y_{\mu\lambda}  = Y_{\lambda\bar \mu}^* =
Y_{\bar \lambda \bar \mu}.
\end{equation*}
\begin{equation*}
Y_{\lambda\mu}  = \sum_k N_{\lambda\mu}^\nu \frac{\omega_\lambda\omega_\mu}
{\omega_\nu} d(\nu) .
\end{equation*}
\end{lemma}
We note that one may take the second equation in the above lemma as the
definition of $Y_{\lambda\mu}$.\par
Define
$a := \sum_i d_{\rho_i}^2 \omega_{\rho_i}^{-1}$.
If the matrix $(Y_{\mu\nu})$ is invertible,
by Proposition on P.351 of \cite{R2} $a$ satisfies
$|a|^2 = \sum_\lambda d(\lambda)^2$.
\begin{definition}\label{c0'}
Let $a= |a| \exp(-2\pi i \frac{c_0}{8})$ where  $c_0\in {\mathbb R}$ and $c_0$
is well defined ${\rm mod} \ 8\mathbb Z$.
\end{definition}
Define matrices
\begin{equation}\label{Smatrix}
S:= |a|^{-1} Y, T:=  C {\rm Diag}(\omega_{\lambda})
\end{equation}
where \[C:= \label{dims} \exp(-2\pi i \frac{c_0}{24}).\]
Then these matrices satisfy (cf. \cite{R2}):
\begin{lemma}\label{Sprop}
\begin{align*}
SS^{\dag} & = TT^{\dag} ={\rm id},  \\
STS &= T^{-1}ST^{-1},  \\
S^2 & =\hat{C},\\
 T\hat{C} & =\hat{C}T,
\end{align*}

where $\hat{C}_{\lambda\mu} = \delta_{\lambda\bar \mu}$
is the conjugation matrix.
\end{lemma}
Moreover
\begin{equation}\label{Verlinde}
N_{\lambda\mu}^\nu = \sum_\delta \frac{S_{\lambda\delta} S_{\mu\delta}
S_{\nu\delta}^*}{S_{1\delta}}. \
\end{equation}
is known as Verlinde formula.
The commutative algebra  generated by $\lambda$'s with structure
constants $N_{\lambda\mu}^\nu$ is called {\bf fusion algebra} of
$\A$. If $Y$ is invertible,
it follows from Lemma \ref{Sprop}, (\ref{Verlinde})  that any nontrivial
irreducible representation
of the fusion algebra is of the form
$\lambda\rightarrow \frac{S_{\lambda\mu}}{S_{1\mu}}$ for some
$\mu$. \par
\subsection{Induced endomorphisms}\label{inductionsection}
Suppose that $\ro\in \End(M)$ has the property that
$\gamma=\rho\bar\rho \in \Delta_\A$. By \S2.7 of \cite{LR}, we can
find two isometries $v_1\in \Hom(\gamma,\gamma^2), w_1\in \Hom
(1,\gamma)$\footnote{We use $v_1,w_1$ instead of $v,w$ here since
$v, w$ are used to denote sectors in Section \ref{typea}.} such that
$\bar\rho(M)$ and $v_1$ generate $M$ and
\begin{align*}
v_1^* w_1 & = v_1^* \gamma(w_1) = d_\rho^{-1} \\
v_1v_1 & = \gamma(v_1) v_1  \\
%\epsilon(\ra,\r) w_1 & = w_1 \label{c}
\end{align*}
By Thm. 4.9 of \cite{LR}, we shall say that $\rho$ is {\it local }
if
\begin{align}\label{local}
v_1^* w_1 & = v_1^* \gamma(w_1) = d_\rho^{-1} \\
v_1v_1 & = \gamma(v_1) v_1  \\
\bar\rho(\epsilon(\gamma,\gamma)) v_1 & = v_1
\end{align}
Note that if $\rho$ is local, then
\begin{equation}\label{local=1}
\om_\mu =1, \forall \mu\prec \ro\bar\ro \end{equation}

For each (not necessarily irreducible)  $\lambda\in \Delta_\A,$ let
$\e{(\lambda,\gamma)}
 :   \lambda\gamma \rightarrow
\gamma\lambda$ (resp.  $\tilde\e{(\lambda,\gamma)}$), be the
positive (resp. negative) braiding operator as defined in Section
1.4 of \cite{Xb}. Denote by $\lambda_\e \in$ End$(M)$ which is
defined by
\begin{align*}
\lambda_\e(x) :&= ad(\e(\lambda,\gamma))\lambda(x)=
\e(\lambda,\gamma) \lambda(x) \e(\lambda,\gamma)^* \\
\lambda_{\tilde\e}(x) :&= ad(\tilde\e(\lambda,\gamma))\lambda(x)=
\tilde\e(\lambda,\gamma)^* \lambda(x) \tilde\e(\lambda,\gamma)^*,
\forall x\in M.
\end{align*}
 By (1) of Theorem 3.1 of \cite{Xb},  $\lambda_{\e} \rho (M)
\subset \rho(M), \lambda_{\tilde\e}  \rho (M) \subset \rho(M),$
hence the following definition makes sense\footnote{We have changed
the notations $a_\la, \tilde a_\la$ of \cite{Xb} to $\tilde a_\la,
a_\la$ of this paper to make some of the formulas such as equation
(\ref{aa2}) simpler.}. \bdef\label{ala} If $\la\in \Delta_\A$ define
two elements of $\End(M)$ by
$$
a^\ro_\lambda(m):= \rho^{-1} (\lambda_\e  \rho (m)) , \ \tilde
a^\ro_\lambda(m):= \rho^{-1} (\lambda_{\tilde\e}  \rho (m)), \forall
m\in M.
$$
$a_\lambda^\ro$ (resp. $\tilde a_\lambda^\ro$) will be referred to
as positive (resp. negative) induction of $\la$ with respect to
$\ro.$

\ede
\begin{remark}
For simplicity we will use $a_\la, \tilde a_\la$ to denote
$a_\lambda^\ro, \tilde a_\lambda^\ro$ when it is clear that
inductions are with respect to the same $\ro.$
\end{remark}

The endomorphisms $a_\lambda$ are called braided endomorphisms in
\cite{Xb}  due to its braiding properties (cf.  (2) of Corollary 3.4
in \cite{Xb}), and enjoy an interesting set of properties (cf.
Section 3 of \cite{Xb}). Though \cite{Xb} focus on the local case
\footnote{ As we will see in Prop. \ref{loclr}, the induction with
respect to non-local $\ro$ is closely related to induction with
respect to certain local $\ro'$ related to $\ro.$} which was clearly
the most interesting case in terms of producing subfactors, as
observed in
 \cite{BE1}, \cite{BE2}, \cite{BE3}, \cite{BE4}
 that many of the arguments in \cite{Xb} can be
generalized. These properties are also studied in a slightly
different context in \cite{BE1}, \cite{BE2}, \cite{BE3}. In these
papers,  the induction is between $M$ and a subfactor $N$ of $M$
,while the induction above is on the same algebra. A dictionary
between our notations here and these papers has been set up in
\cite{Xj} which simply use an isomorphism between $N$ and $M$. Here
one has a choice to use this isomorphism to translate all
endomorphisms of $N$ to endomorphims of $M$, or equivalently all
endomorphims of $M$ to endomorphims of $N$. In \cite{Xj} the later
choice is made (Hence $M$ in \cite{Xj} will be our $N$ below). Here
we make the first choice which makes the dictionary slightly
simpler. Our dictionary here is equivalent to that of \cite{Xj}. Set
$N=\bar\rho(M)$. In the following the notations from \cite{BE1} will
be given a subscript BE. The formulas are :
%Denote by $\lambda_\sigma \in$ End$(M)$ which is defined by
%$\lambda_\sigma(x) = \sigma^* \lambda(x) \sigma$ for any $x \in M$.
%$\gamma_{BE}$ (resp. $\theta_{BE}$) of [5] corresponds to
%our $\alpha$ (resp. $\gamma$ )
%,  $\lambda_{BE}\in End (
%\pi^0(N(I_0)))$ in [5] corresponds to our
%$\phi_{I_0} \lambda \phi_{I_0}^{-1}$, and
%$$\epsilon (\lambda_{BE}, \theta_{BE})$
%of [5] corresponds to our $\phi_{I_0}(\e{\lambda\gamma})$. Also
%note that
%$\rho^{-1} \phi_{I_0^c} (\epsilon (\lambda_{BE}, \theta_{BE}))
%= \rho^{-1} \gamma ( \sigma_{\lambda\gamma}) = \bar\rho (\sigma_
%{\lambda\gamma})$ by (2.1.3) and $\rho \bar\rho= \gamma$.
%We summarize this by the following formulas:
\begin{align}\label{aa1}
\rho\res N=i_{BE}, \ & \bar\ro\ro\res N=\bar i_{BE}i_{BE}, \\
\gamma = \bar\rho^{-1} \theta_{BE} \bar\rho, \ &
\bar\ro \ro= \gamma_{BE}, \\
\lambda = \bar\ro^{-1} \lambda_{BE} \bar\ro, \ \ & \e(\lambda, \mu)
= \bar\rho (\e^{+}(\lambda_{BE}, \mu_{BE})) \\
& \tilde\e(\lambda, \mu) = \bar\rho (\e^{-}(\lambda_{BE}, \mu_{BE}))
\end{align}
The dictionary between $a_\lambda \in End (M)$ in definition
\ref{ala} and $\alpha_\lambda^{-}$ as in Definition 3.3, Definition
3.5 of \cite{BE1}  are given by:
\begin{equation}\label{aa2}
a_\lambda  =\alpha^+_{\lambda_{BE}}, \tilde a_\lambda
=\alpha^{-}_{\lambda_{BE}}
\end{equation}
The above formulas will be referred to as our {\it dictionary}
between the notations of \cite{Xb} and that of \cite{BE1}. The proof
is the same as that of \cite{Xj}. Using this dictionary one can
easily translate results of \cite{Xb} into the settings of
\cite{BE1},  \cite{BE2}, \cite{BE3}, \cite{BE4},\cite{BEK1},
\cite{BEK2}and vice versa. First we summarize a few properties from
\cite{Xb} which will be used in this paper: (cf. Th. 3.1 , Co. 3.2
and Th. 3.3 of \cite{Xb} ): \bprop\label{xua} (1). The maps
$[\lambda] \rightarrow [a_\lambda], [\lambda]\rightarrow [\tilde
a_\lambda]$ are ring homomorphisms;\par (2) $ a_\lambda \bar \rho
=\tilde a_\lambda \bar \rho =\bar \rho \lambda$;\par (3) When
$\ro\bar\ro$ is local, $ \l a_\la,  a_\mu \ra = \l \tilde a_\la,
\tilde a_\mu \ra = \l a_\la\bar \ro, a_\mu\bar\ro \ra =\l \tilde
a_\la\bar \ro, \tilde a_\mu\bar\ro \ra ;$\par (4) (3) remains valid
if $a_\la, a_\mu$ (resp. $\tilde a_\la, \tilde a_\mu$) are replaced
by their subsectors.

\eprop

\bdef\label{hro} $H_\ro$ is a finite dimensional vector space over
$\Co$ with orthonormal basis consisting of irreducible sectors of
$[\lambda\rho],\forall \lambda\in \Delta_\A.$ \ede
\par

$[\lambda]$ acts linearly on $H_\ro$ by $[\lambda] [a]= \sum_b \l
\lambda a,b\ra [b]$ where $[b]$ are elements in the basis of
$H_\ro.$ \footnote{By abuse of notation, in this paper we use
$\sum_b$ to denote the sum over the basis $[b]$ in $H_\ro$.}  By
abuse of notation, we use $[\lambda]$ to denote the corresponding
matrix relative to the basis of $H_\ro.$ By definition these
matrices are normal and commuting, so they can be simultaneously
diagonalized. Recall the irreducible representations of the fusion
algebra of $\A$ are given by
$$
\lambda \rightarrow \frac{S_{\lambda \mu}}{S_{1\mu}}.
$$
\bdef\label{orthogonalphig} Assume $\l \la a,b\ra=
\sum_{{\mu,i}\in(\text{\rm Exp})} \frac{S_{\lambda \mu}}{S_{1\mu}}
\cdot \phi_a^{(\mu,i)} \phi_b^{(\mu,i)^*}$ where $ \phi_a^{(\mu,i)}$
are normalized orthogonal eigenvectors of $[\lambda]$  with
eigenvalue $\frac{S_{\lambda \mu}}{S_{1\mu}},$ $\Exp$ is a set of
$\mu,i$'s
 and $i$ is an
index indicating the multiplicity of  $\mu$. Recall if a
representation is denoted by $1$, it will always be the vacuum
representation. \ede

The following lemma is elementary:

\blem\label{ephi} (1):
\[\sum_b d_b^2=\frac{1}{S_{11}^2}\]
where the sum is over the basis of $H_\ro$. The vacuum appears once
in $\Exp$ and
\[
\phi_a^{(1)}= S_{11} d_a; \]\par (2) \[\sum_i
\frac{\phi_a^{(\la,i)}{{\phi_a^{(\la,i)}}^*}}{S_{1\la}^2}= \sum_\nu
\l \bar\nu a,b\ra \frac{S_{\nu\la}}{S_{1\la}}\] where if $\la$ does
not appear in $\Exp$ then the righthand side is zero.

\elem

\prf Ad (1): By definition we have
\[[a\bar\rho]= \sum_{\la}\l a\bar\rho,\la\ra[\la]=\sum_{\la}\l a,\la\rho\ra[\la]
\] where in the second $=$ we have used Frobenius reciprocity. Hence
\[d_a d_{\bar\ro}= \sum_\lambda \l a\bar\rho,\la\ra d_\la\] and we
obtain
\[
\sum_\la d_\la^2 = \sum_{\la,a} \l a\bar\ro,\la\ra d_\la d_a/d_\ro =
\sum_a d_a^2
\]
(2) follows from definition and orthogality of $S$ matrix. \qed

\par \subsection{ Relative braidings}\label{relbraidsection}
In \cite{Xb}, commutativity among subsectors of $a_\la,\tilde a_\mu$
were studied.  We record these results in the following for later
use: \blem\label{acommut} (1) Let $[b]$ (resp. $[b']$) be any
subsector of $a_\la$ (resp. $\tilde a_\la$). Then
\[
[a_\mu b]=[ba_\mu], [\tilde a_\mu b']=[b'\tilde a_\mu]\forall \mu,
[bb']=[bb'];
\]\par
(2) Let $[b]$ be a subsector of $a_\mu\tilde a_\la,$ then $[a_\nu
b]=[b a_\nu],[\tilde a_\nu b]=[b \tilde a_\nu], \forall \nu.$ \elem

\prf (1) follows from (1) of Th. 3.6 and Lemma 3.3 of \cite{Xb}. (2)
follows from the proof of Lemma 3.3 of \cite{Xb}. Also cf. Lemma
3.20 of \cite{BE3}.\qed

In the proof of these commutativity relations in \cite{Xb}, an
implicit use of relative braidings was made. These braidings are
further studied in \cite{BE2}, \cite{BE3} and let us recall their
properties in our setting by using our dictionary (\ref{aa1}),
(\ref{aa2}).

Let $\tilde \beta, \delta\in \text{\rm End}(M)$ be subsectors of
$\tilde a_\lambda$ and $a_\mu$. By Lemma 3.3 of \cite{Xb}, $[\tilde
\beta]$ and $[\delta]$ commute. Denote by $\epsilon_r (\tilde\beta,
\delta)$ given by:
\begin{align}\label{eq6}
\epsilon_r (\tilde\beta, \delta):&= s^* a_\mu (t^*)
\bar\rho(\sigma_{\lambda\mu})
\tilde a_\lambda (s) t \in \text{\rm Hom} (\beta \delta, \delta\beta) \\
\epsilon_r (\delta, \tilde\beta):&= \epsilon_r (\tilde\beta,
\delta)^{-1},
\end{align}
 with isometries $t\in  \text{\rm Hom} (\tilde \beta, \tilde
a_\lambda)$ and $s\in \text{\rm Hom}  (\delta,  a_\mu)$. Also
$$
\epsilon_r (\tilde a_\lambda, a_\mu) = \bar\rho(\sigma_{\lambda\mu})
\ , \ \epsilon_r ( a_\lambda, \tilde a_\mu) = \bar\rho(\tilde
\sigma_{\lambda\mu})
$$
\par
\blem\label{rb1}
The operator $\epsilon_r (\beta, \delta)$ defined
above does not depend on $\lambda, \mu$ and the isometries $s,t$ in
the sense that, if there are isometries $x\in Hom (\beta, \tilde
a_\nu)$ and $y\in \text{\rm Hom}( \delta, a_{\delta_1})$, then
$$
\epsilon_r (\beta, \delta) =  s^* a_{\delta_1}
 (t^*) \bar\rho(\sigma_{\nu\lambda_1})
\tilde a_\nu (y) x
$$
\elem
\par
\blem\label{relb} The system of unitaries of Eq. (\ref{eq6})
provides a relative braiding between representative endomorphisms of
subsectors of $\tilde a_\lambda$ and $a_\mu$ in the sense that , if
$\beta, \delta, \omega, \xi$ are subsectors of $[\tilde a_\lambda],
[ a_\mu], [\tilde a_\nu], [a_{\delta_1}]$, respectively, then we
have initial conditions
$$
\epsilon_r (\text{\rm id}_M, \delta) = \epsilon_r (\beta, \text{\rm
id}_M) = 1 ,$$ compositions rules
$$
\epsilon_r ( \beta \omega, \delta) = \epsilon_r ( \beta,\delta)
\beta( \epsilon_r (\omega, \delta)), \epsilon_r ( \beta,  \delta
\xi) = \delta (\epsilon_r ( \beta, \xi))
 \epsilon_r (\beta, \delta)
,$$ and naturality
$$
\delta (q_+) \epsilon_r (\beta, \delta) =\epsilon_r (\omega, \delta)
q_+, q_-,  \epsilon_r (\beta, \delta) = \epsilon_r (\beta, \xi)
\beta (q_-)
$$ whenever
$q_+ \in \text{\rm Hom}(\beta, \omega)$ and $q_-\in Hom (\delta,
\xi)$. \par For the collection of $\beta,\delta$ such that
$\beta\prec a_\la, \beta\prec \tilde a_\la$ and $\delta\prec a_\mu,
\delta\prec \tilde a_\mu$ for some (varying) $\la,\mu \in
\Delta_\a,$ the unitaries $\e_r(\beta,\delta), \e_r(\delta,\beta)$
defines a braiding: i.e., they verify YBE and BFE in Prop.
\ref{bfe}.

\elem \blem\label{313} Let $r\in \Hom(\lambda_3,
\lambda_1\lambda_2).$ Then
$$\bar\rho(r)\in \Hom(a_{\la_3},
a_{\la_1}a_{\la_2})\cap \Hom(\tilde a_{\la_3}, \tilde
a_{\la_1}\tilde a_{\la_2}).
$$ \elem \prf When $\rho\bar\ro$ is local,
the lemma follows from Th. 3.3 of \cite{Xb}. Let us prove the
general case. Since $a_\lambda\bar\ro= \bar\ro \lambda$, we have
$\bar\rho(r)\in \Hom(a_{\la_3}\bar\rho, a_{\la_1\la_2}\bar\ro).$
Since $M$ is generated by $\bar\ro(M), v_1$, to finish the proof we
just need to check that
\[
\bar\ro(r)a_{\lambda_3}(v_1)=a_{\la_1\la_2}(v_1)\bar\rho(r)
\]
Since $\rho$ is one to one, applying $\ro$ to the above equation it
is sufficient to check that
\[\gamma(r)\ro a_{\lambda_3}(v_1)=\ro a_{\la_1\la_2}(v_1)\gamma(r)
\]
Using $\rho a_\la= \e(\la,\gamma) \la\ro \e(\la,\gamma)^*$, one can
check directly that this equation follows from  Prop. \ref{bfe}.
\qed

The following is Lemma 3.25 of \cite{BE1} in our setting:
\blem\label{325} If $r\in \Hom (\bar\ro\la,\bar\ro\mu),$ then
\[
r\bar\ro(\e(\mu_1,\la))= \bar\ro(\e(\mu_1,\la))a_{\mu_1}(r),
r\bar\ro(\tilde \e(\mu_1,\la))= \bar\ro(\tilde \e(\mu_1,\la))\tilde
a_{\mu_1}(r).
\]
\elem Following \cite{BEK1} we define

\bdef\label{zmatrix} For $\la,\mu\in \Delta_\A$, $Z_{\la\mu}:= \l
a_\lambda, \tilde a_\mu\ra.$ \ede
 We can now translate Th. 5.7 and Th. 6.12 of \cite{BEK1} into
our setting: \bprop\label{spec} (1)$\mu$ appears in $\Exp$ as
defined in Definition \ref{orthogonalphig} with multiplicity
$Z_{\mu\mu};$ \par (2) $Z_{\la\mu}$ as a matrix commute with $S,T$
matrices as defined in equation (\ref{Smatrix}).

\eprop By Lemma \ref{ephi} and Prop. \ref{spec} we have the
following: \blem\label{znot0} If
\[\sum_\nu\l\bar\nu a,b\ra\frac{S_{\nu\la}}{S_{1\la}}\ne 0,\]
then $\l a_\lambda, \tilde a_\lambda\ra \geq 1 $ \elem

The following follows from Prop. 3.1 of \cite{BEK1}:

\blem\label{ibfe} For any $\la\in \Delta_\A, b \in H_\ro$ we have
$\e(\la,b\bar\ro) \in \Hom(\la b,b a_\la), \tilde\e(\la,b\bar\ro)
\in \Hom(\la b,b \tilde a_\la).$\elem

Later we will consider the following analogue of $S$-matrix using
relative braidings. Suppose that $T_\mu\in \Hom(a_\mu,\tilde a_\mu),
\forall \mu\in \Delta_\A$ ($T_\mu$ can be zero).
\bdef\label{generalphi} For $\mu\in \Delta_\A, b\in H_\ro$
irreducible , define
\[ \psi_b^{(T_\mu)}:=  S_{11} d_bd_\mu \phi_\mu(
\e(b\bar\ro,\mu)b(T_\mu)\e(\mu,b\bar\ro)).
\]\ede
\blem\label{epsi} (1): $\psi_b^{(T_\mu)}$ depends only on $[b];$
\par
(2) \[ \sum_b {\psi_b^{(T_\mu)}}^* [b]
\] is either zero or an eigenvector of $[\lambda]$ with eigenvalue
$\frac{S_{\la\mu}}{S_{1\mu}},$ and $\sum_b \psi_b^{(T_\mu)} d_b =0$
unless $[\mu]=[1];$ \par (3) If $T_\mu, T_{\bar\mu}$ are unitaries,
and for any irreducible $\la\prec \mu\bar \mu,$ $1\prec a_\la$ iff
$[\la]=[1],$ then $|\sum_b
\psi_b^{(T_\mu)}\psi_b^{(T_{\bar\mu})}|=1;$ \par (4) If  $T_\mu$ is
unitary then $|\psi_b^{(T_\mu)}|\leq S_{11} d_\mu d_b.$

\elem \prf Ad(1): Suppose that $[b_1]=[b]$ and let $U\in
\Hom(b_1,b)$ be a unitary. We have
\begin{align*}
\psi_b^{(T_\mu)}&=  S_{11} d_bd_\mu \phi_\mu(
\e(b\bar\ro,\mu)b(T_\mu)\e(\mu,b\bar\ro)) \\
&=  S_{11} d_bd_\mu \phi_\mu(\mu(U^*)
\e(b\bar\ro,\mu)b\bar\ro(T_\mu)\e(\mu,b\bar\ro)\mu(U)) \\
&=S_{11} d_bd_\mu \phi_\mu(
\e(b_1\bar\ro,\mu)U^*b(T_\mu)U\e(\mu,b_1\bar\ro)) \\
&=S_{11} d_bd_\mu \phi_\mu(
\e(b_1\bar\ro,\mu)b_1(T_\mu)\e(\mu,b_1\bar\ro)) \\
&=\psi_{b_1}^{(T_\mu)}
\end{align*}
Where we have used   BFE of Prop. \ref{bfe} in the third =.
\par Ad (2): Let $t_{b,i} \in \Hom(b,\bar\la b')$ be isometries such
that $\sum_i t_{b,i} t_{b,i}^* =1.$ Then
\[\sum_b \psi_b^{(T_\mu)} \l b,\bar\la b'\ra=
\sum_{b,i} S_{11} d_\mu d_\la d_{b'}\phi_{\bar\la}\phi_\mu
(\mu(t_{b,i}) \e(b\bar\ro,\mu)^*
b(T_\mu)\e(\mu,b\bar\ro)\mu(t_{b,i}^*))
\] where we have used equation
(\ref{phiprop}). By Prop. \ref{bfe} we have
\begin{align*}
&\sum_{b,i} S_{11} d_\mu d_\la d_{b'} \phi_{\bar\la}
\phi_\mu(\mu(t_{b,i}) \e(b\bar\ro,\mu)
b(T_\mu)\e(\mu,b\bar\ro)\mu(t_{b,i}^*)) \\
& = S_{11} d_\mu d_\la d_{b'} \phi_{\bar\la} \phi_\mu(
\e(\bar\la b'\bar\ro,\mu) b(T_\mu)\e(\mu,\bar\la b'\bar\ro))\\
&= \frac{S_{\bar\la \mu}}{S_{1\mu}} \psi_{b'}^{(T_\mu)}
\end{align*}
 Hence
\[
\sum_b [\la] {\psi_b^{(T_\mu)}}^* [b]= \sum_{b,b'}
{\psi_b^{(T_\mu)}}^*\l b,\bar\la b'\ra [b']=
\frac{S_{\la\mu}}{S_{1\mu}}\sum_{b'} {\psi_{b'}^{(T_\mu)}}^* [b'].
\]
By (1) of Lemma \ref{ephi}  we conclude that $\sum_b
\psi_b^{(T_\mu)} d_b =0$ unless $[\mu]=[1].$ \par

Ad (3): Let $t_{\la,i}\in \Hom(\la,\mu\bar\mu)$ be isometries such
that $\sum_{\la,i} t_{\la,i} t_{\la,i}^* =1.$ Then
\begin{align*}
 \psi_b^{(T_\mu)}\psi_b^{(T_{\bar\mu})}&= S_{11} d_bd_\mu \phi_{\bar\mu}(\psi_b^{(T_\mu)}
\e(b\bar\ro,\bar\mu)b(T_\mu)\e(\bar\mu,b\bar\ro)) \\
&=S_{11}^2 d_b^2 d_\mu \phi_{\mu\bar\mu}(
\e(b\bar\ro,\mu\bar\mu)b(T_\mu
a_\mu(T_{\bar\mu}))\e(\mu\bar\mu,b\bar\ro))\\
&= S_{11}^2 d_b \sum_{\la,i} d_bd_\la\phi_\la (\e(b\bar\ro, \la) b
(\bar\ro(t_{\la,i})^* T_\mu
a_\mu(T_{\bar\mu})\bar\ro(t_{\la,i}))\e(\la,b\bar\ro))
\end{align*}
Where we have used equation (\ref{phiprop}) and Lemma \ref{ibfe} in
the second = and BFE of Prop. \ref{bfe} in the third =. By (2) of
Lemma \ref{epsi}
\[\sum_b d_b d_bd_\la\phi_\la (\e(b\bar\ro, \la)
b (\bar\ro(t_{\la,i})^* T_\mu
a_\mu(T_{\bar\mu})\bar\ro(t_{\la,i}))\e(\la,b\bar\ro)) =0\] unless
$[\la]=[1].$ Denote by $t_1\in \Hom (1,\mu\bar\mu)$ the unique (up
to scalar) isometry. Then we have (Recall we always identify the
center of $M$ with $\Co$)
 \[
\sum_b \psi_b^{(T_\mu)}\psi_b^{(T_{\bar\mu})}=  \bar\ro(t_{1})^*
T_\mu a_\mu(T_{\bar\mu})\bar\ro(t_{1})
\]
On the other hand since $T_\mu, T_{\bar \mu}$ are unitaries, we have
\[
\sum_{\la,i} \bar\ro(t_{1})^*T_\mu
a_\mu(T_{\bar\mu})\bar\ro(t_{\la,i})\bar\ro(t_{\la,i})^*
a_\mu(T_{\bar\mu})^* T_\mu^*\bar\ro(t_{1}) =1
\]
Since $\bar\ro(t_{1})^*T_\mu a_\mu(T_{\bar\mu})\bar\ro(t_{\la,i})\in
\Hom( a_\la, 1),$ by assumption it is $0$ unless $[\la]=[1].$ We
conclude that $|\bar\ro(t_{1})^* T_\mu
a_\mu(T_{\bar\mu})\bar\ro(t_{1})|=1$ and (3) is proved. (4) follows
since $\phi_\mu$ is completely positive.

\qed

Using equation (\ref{aa1}), (\ref{aa2}), the following is a
translation of Prop. 3.2 and Th. 4.7 of \cite{BE4} into our setting:
\bprop\label{loclr} Suppose that $\ro\bar\ro \in \Delta$. Then:\par
(1) $\ro$ is local iff $\l 1, a_\mu\ra =\l \ro \bar\ro, \mu\ra,
\forall \mu\in \Delta_\A;$ \par (2) \[ \ro= \ro'\ro''=
\tilde\ro'\tilde\ro'' \] where $\ro',\ro'', \tilde\ro', \tilde\ro''
\in \End (M)$, and $\ro',\tilde\ro'$ are local which verifies
\begin{align*}
\l \ro'\bar\ro',\mu\ra &= \l 1,a_\mu\ra =\l 1, a_\mu^{\ro'}\ra
\\
\l \tilde\ro'\overline{\tilde\ro'},\mu\ra &= \l 1,\tilde a_\mu\ra
=\l 1, \tilde a_\mu^{\tilde \ro'}\ra
\end{align*}
$\forall \mu\in \Delta_\A.$ \eprop The following Lemma is Prop. 3.23
of \cite{BE1} (The proof was also implicitly contained in the proof
of Lemma 3.2 of \cite{Xb}):

\blem\label{loc2} If $\rho\bar\ro$ is local, then $[a_\la]=[\tilde
a_\la]$ iff $\e(\la,\rho\bar\rho)\e(\ro\bar\ro,\la) =1.$ \elem

%\subsection{Type A}
\subsection {Jones-Wassermann subfactors from representation of Loop
groups}\label{typea} Let $G= SU(n)$. We denote $LG$ the group of
smooth maps $f: S^1 \mapsto G$ under pointwise multiplication. The
diffeomorphism group of the circle $\text{\rm Diff} S^1 $ is
naturally a subgroup of $\text{\rm Aut}(LG)$ with the action given
by reparametrization. In particular the group of rotations
$\text{\rm Rot}S^1 \simeq U(1)$ acts on $LG$. We will be interested
in the projective unitary representation $\pi : LG \rightarrow U(H)$
that are both irreducible and have positive energy. This means that
$\pi $ should extend to $LG\rtimes \text{\rm Rot}\ S^1$ so that
$H=\oplus _{n\geq 0} H(n)$, where the $H(n)$ are the eigenspace for
the action of $\text{\rm Rot}S^1$, i.e., $r_\theta \xi = \exp(i n
\theta)$ for $\theta \in H(n)$ and $\text{\rm dim}\ H(n) < \infty $
with $H(0) \neq 0$. It follows from \cite{PS} that for fixed level
$k$ which is a positive integer, there are only finite number of
such irreducible representations indexed by the finite set
$$
 P_{++}^{k}
= \bigg \{ \lambda \in P \mid \lambda = \sum _{i=1, \cdots , n-1}
\lambda _i \Lambda _i , \lambda _i \geq 0\, , \sum _{i=1, \cdots ,
n-1} \lambda _i \leq k \bigg \}
$$
where $P$ is the weight lattice of $SU(n)$ and $\Lambda _i$ are the
fundamental weights. We will write $\la=(\la_1,...,\la_{n-1}),
\la_0= k-\sum_{1\leq i\leq n-1} \la_i$ and refer to
$\la_0,...,\la_{n-1}$ as components of $\la.$

We will use $\Lambda_0$ or simply $1$  to denote the trivial
representation of $SU(n)$. For $\lambda , \mu , \nu \in P_{++}^{k}$,
define $N_{\lambda \mu}^\nu  = \sum _{\delta \in P_{++}^{k}
}S_\lambda ^{(\delta)} S_\mu ^{(\delta)} S_\nu
^{(\delta*)}/S_{\Lambda_0}^{(\delta})$ where $S_\lambda ^{(\delta)}$
is given by the Kac-Peterson formula (cf. equation (\ref{kacp})
below for an equivalent formula):
$$
S_\lambda ^{(\delta)} = c \sum _{w\in S_n} \varepsilon _w \exp
(iw(\delta) \cdot \lambda 2 \pi /n)
$$
where $\varepsilon _w = \text{\rm det}(w)$ and $c$ is a
normalization constant fixed by the requirement that
$S_\mu^{(\delta)}$ is an orthonormal system. It is shown in
\cite{Kac2} P. 288 that $N_{\lambda \mu}^\nu $ are non-negative
integers. Moreover, define $ Gr(C_k)$ to be the ring whose basis are
elements of $ P_{++}^{k}$ with structure constants $N_{\lambda
\mu}^\nu $.
  The natural involution $*$ on $ P_{++}^{k}$ is
defined by $\lambda \mapsto \lambda ^* =$ the conjugate of $\lambda
$ as representation of $SU(n)$.\par

We shall also denote $S_{\Lambda _0}^{(\Lambda)}$ by $S_1^{(\Lambda
)}$. Define $d_\lambda = \frac {S_1^{(\lambda )}}{S_1^{(\Lambda
_0)}}$. We shall call $(S_\nu ^{(\delta )})$ the $S$-matrix of
$LSU(n)$ at level $k$. \par
 We shall encounter the $\Bbb Z_n$
group of automorphisms of this set of weights, generated by
$$
\sigma : \lambda = (\lambda_1, \lambda_2, \cdots , \lambda_{n-1})
\rightarrow \sigma(\lambda) = ( k -1- \lambda_1 -\cdots
\lambda_{n-1}, \lambda_1, \cdots , \lambda_{n-2}).
$$
Define  $\col(\lambda) = \Sigma_i (\lambda_i - 1) i $. The central
element  $ \exp \frac{2\pi i}{n}$ of $SU(n)$ acts on representation
of $SU(n)$ labeled by $\lambda$ as $\exp( \frac{2\pi i
\col(\lambda)}{n})$. The irreducible positive energy representations
of $ L SU(n)$ at level $k$ give rise to an irreducible conformal net
$\A$ (cf. \cite{KLX}) and its covariant representations. We will use
$\la=(\la_1,...\la_{n-1})$ to denote irreducible representations of
$\A$ and also the corresponding endomorphism of $M=\A(I).$

All the sectors $[\lambda]$ with $\lambda$ irreducible generate the
fusion ring of $\A.$
\par For $\lambda$ irreducible, the univalence $\om_\lambda$ is given
by an explicit formula (cf. 9.4 of [PS]). Let us first define
$h_\lambda = \frac {c_2(\lambda)}{k+n}$ where $c_2(\lambda)$ is the
value of Casimir operator on representation of $SU(n)$ labeled by
dominant weight $\lambda$.
 $h_\lambda$ is usually called the conformal dimension. Then
we have: $\om_\lambda = exp({2\pi i} h_\lambda)$. The conformal
dimension of $\lambda=(\la_1,...,\la_{n-1})$ is given by
\begin{equation}\label{cdim} h_\lambda= \frac{1}{2n(k+n)}\sum_{1\leq
i\leq n-1} i(n-i) \la_i^2 + \frac{1}{n(k+n)}\sum_{1\leq j\leq i\leq
n-1}j (n-i)\la_j\la_i + \frac{1}{2(k+n)}\sum_{1\leq j\leq n-1}
j(n-j) \la_j \end{equation} The following form of Kac-Peterson
formula for $S$ matrix will be used later:
\begin{equation}\label{kacp}
\frac{S_{\la\mu}}{S_{1\mu}}= Ch_{\la'} (x_1,...,x_{n-1}, 1)
\end{equation}
Where $ch_{\la'}$ is the character associated with finite
irreducible  representation of $SU(n)$ labeled by $\la,$ and
$x_i=\exp(-2\pi i\frac{\mu_i'}{k+n}), \mu_i'=\sum_ {i\leq j\leq n-1}
(\mu_j+1), 1\leq i\leq n-1.$ \par

It follows that $S$ matrix verifies:

\begin{equation}\label{ssy}
S_{\la\om^i(\mu)}= \exp(\frac{2\pi i \col(\la)}{n}) S_{\la\mu}
\end{equation}
\par The
following  result is proved in \cite{Wass} (See Corollary 1 of
Chapter V in \cite{Wass}).

\bthm\label{wass}  Each $\lambda \in  P_{++}^{(k)}$ has finite index
with index value $d_\lambda ^2$.  The fusion ring generated by all
$\lambda \in P_{++}^{(k)}$ is isomorphic to $ Gr(C_k)$. \ethm

\begin{remark}\label{jxx}
The subfactors in the above theorem are called Jones-Wassermann
subfactors after the authors who first studied them (cf.
\cite{Jh1},\cite{Wass}).
\end{remark}
\begin{definition}
$v:=(1,0,...,0), \ad:=(1,0,...,0,1), \om^i =k \Lambda_i, 0\leq i\leq
n-1. $ \ede The following is observed in \cite{GW}:
\blem\label{fusionrule} Let $(0,...,0,1,0,...0)$ be the $i$-th
($1\leq i\leq n-1$) fundamental weight. Then
$[(0,...,0,1,0,...0)\la]$ are determined as follows:
$\mu\prec(0,...,0,1,0,...0)\la$ iff when the Young diagram of $\mu$
can be obtained from Young diagram of $\la$ by adding $i$ boxes on
$i$ different rows of $\la$, and such $\mu$ appears in
$[(0,...,0,1,0,...0)\la]$ only once. \elem

\blem\label{fusionv} (1) If $[\la]\neq \omega^i$ for some $0\leq
i\leq n-1,$ then $\ad\prec \la\bar\la;$\par (2) If
$\lambda_1\lambda_2$ is irreducible, then either $\lambda_1$ or
$\lambda_2=\omega^i$ for some $0\leq i\leq n-1$. \elem

\prf By Lemma \ref{fusionrule} we have that
\[\l v \la,v\la \ra=1\] iff $\la =\omega^i$ for some $0\leq i\leq
n-1$. By Frobenius reciprocity
\[\l v \la,v\la \ra= \l 1+ \ad,\la\bar\la\ra= 1+ \l \ad,\la\bar\la\ra\]
Hence \[\l \ad,\la\bar\la\ra=0\] iff $\la =\omega^i$ for some $0\leq
i\leq n-1$. If $\lambda_1\lambda_2$ is irreducible, then by
Frobenius reciprocity again we have
\[ \l \la_1\bar\la_1, \la_2\bar\la_2\ra
=1\geq 1+ \l \ad,\la_1\bar\la_1\ra \l \ad,\la_2\bar\la_2\ra
\]
Hence either \[ \l \ad,\la_1\bar\la_1\ra=0\] or
\[ \l\ad,\la_2\bar\la_2\ra=0\] and the lemma follows.
\qed

\blem\label{scalarv} Suppose $\la\in \Delta_\A$ and $\la$ is not
necessarily irreducible. Then
\[\e(\la,v)\e(v,\la) \in \Co \]
iff $[\la]=\sum_{j} [\omega^j]$ where the summation is over a finite
set. \elem

\prf By Prop. \ref{bfe} we have that
\[ \e(v^m,\la)\e(\la,v^m) \in \Co\] for all $m\geq 0$. Since any
irreducible $\mu$ is a subsector of $v^m$ for some $m\geq 0$, by
Lemma \ref{escalar} we have that $\e(\mu,\la_1)\e(\la_1,\mu)\in \Co,
\forall \mu,\la_1\prec \la$. By definition of $S$ matrix we have
$|S_{\mu\la_1}|^2= |S_{1\la_1} d_\mu|^2.$ Sum over $\mu$ we have
$d_{\la_1} = 1,$ i.e., $\la_1$ is an automorphism, and this implies
that $v\la_1$ is irreducible. The lemma now follows from Lemma
\ref{fusionv}. \qed

\blem\label{denseb} For any $m\geq 1, \Hom(v^m,v^m)$ is generated as
an algebra by $1, v^i (\e(v,v)), 1\leq i\leq m-1$. \elem

\prf This is (3) of Lemma 3.1.1 in \cite{Xj} and is essentially
contained in \cite{We}. \qed

Now let $\rho\bar\ro\in \Delta_\A$ where $\A$ is the conformal net
associated with $SU(n)$ at level $k,$ and consider induction with
respect to $\ro$ as defined in Definition \ref{ala}. We have
\blem\label{irredad} (1) $a_v,\tilde a_v$ are always irreducible;
\par (2) $d_{\ad}=1$ iff $k=n=2;$\par (3) If $k\neq n\pm 2,n,$ then
$a_{\ad}, \tilde a_{\ad}$ are irreducible. \elem \prf It is enough
to prove the Lemma for positive induction. The negative induction
case is similar. Assume that $\ro=\ro'\ro''$ as in Prop.
\ref{loclr}, since $\l a_\la, 1\ra= \l \rho'\bar\ro',\la\ra = \l
a_\la^{\rho'}, 1\ra, \forall \la,$ it is enough to prove the Lemma
for induction with respect to $\ro'.$ Hence we may assume that $\ro$
is local. By (3) of Prop. \ref{xua} we have
\[
\l a_v,a_v\ra = \l \ro\bar\ro, v\bar v\ra =1+ \l \ro\bar\ro, \ad\ra
\]
Since $\om_\ad= \exp(\frac{2\pi in}{k+n})\neq 1,$ by equation
(\ref{local=1}) we conclude that $\l \ro\bar\ro, \ad\ra=0$ and (1)
is proved. (2) follows from equation (\ref{kacp}).\par Ad (3) By
Lemma \ref{fusionrule} we  have \[ [\ad^2]= [1]+ 2[\ad] +
[{(2,0,...,0,2)}] + [{(0,1,0,...,1,0)}]+ [{(0,1,0,...,0,2)}] +
[{(2,0,...,0,1,0)}]
\]
By computing the conformal dimensions of the descendants of $\ad^2$
using equation (\ref{cdim}) we have
\[
h_{(2,0,...,0,2)}= \frac{2+2n}{k+n},
h_{(0,1,...,0,2)}=h_{(2,0,...,1,0)} =\frac{2n}{k+n},
h_{(0,1,...,1,0)}= \frac{2n-2}{k+n} \] By equation (\ref{local=1})
we conclude that if $k\neq n\pm2, n,$ then $\l \ad^2, \ro\bar\ro\ra
=1$ and (3) is proved.

\par

\subsection{Induced subfactors from simple current
extensions}\label{orb}
 In this section  we assume that the level
$k=n'n$ where $n'\geq 3$, and $n'$ is an even integer if $n$ is
even. This last condition comes from \cite{Xo}.  For such level it
is shown in \S3 of \cite{BE2} that there is $\rho_o\in \End(M)$ such
that $[\rho_o\bar\rho_o]= \sum_{0\leq i\leq n-1} [\omega^i]$ and
$\ro_o\bar\ro_o$ is local. It also follows from definitions that one
can choose  $\bar\ro_o\rho_o= \sum_{0\leq i\leq n-1} [g^i]$ where
$[g^n]=[1]$ and $[\tilde a_v]=[a_v g]$ (cf. \S6.1 of \cite{KLX}).
Also note that $[a_{\omega^i}]= [1], \forall i.$ The following is a
consequence of Lemma \ref{ephi} and Prop. \ref{xua}:

\blem\label{orthophio} There exists an orthonormal basis $\sum_{a}
\phi_a^{\mu} [a]$ where $\col(\mu)= 0\mod \  n$ and  the sum is over
all irreducible subsectors of $a_\la,\forall \la$ such that
\[
\l a_\la a,b\ra= \sum_{\mu,i,\col(\mu)=0 \mod \  n}
\frac{S_{\la\mu}}{S_{1\mu}} \phi_a^{(\mu,i)} {\phi_b^{(\mu,i)}}^*
\]
\elem The following follows from Cor. 4.9 of \cite{KLX}:
\blem\label{klx} (1) Let $\la$ be irreducible and suppose $l$ is the
smallest positive integer with $[\omega^{l}\la]=[\la]$. Then
$[a_\la]=\sum_{1\leq i\leq l'} [x_i]$ where $l'l=n$ and
$[g^{i}x_1g^{-i}]=[x_i], 1\leq i\leq l', [x_i]\neq [x_j]$ if $i\neq
j.$ Similar statements hold true for $\tilde a_\la;$
\par (2) $\l a_\la,a_\mu\ra \neq 0$ iff $[\la]=[\om^j(\mu)]$ for
some $1\leq j\leq n$ iff $[a_\la]=[a_\mu].$ Similar statements hold
true for $\tilde a_\la, \tilde a_\mu.$

 \elem
Later we will use the following analogue of Lemma \ref{scalarv}:

\blem\label{scalarad} If $\e(\ad, \la)\e(\la,\ad)\in \Co$, then
$[\la]=\sum_j \omega^j$ where the sum is over a finite set of
positive integers.\elem \prf By Prop. \ref{bfe}  and Lemma
\ref{escalar} we have that $\e(\ad^m,\la_1)\e(\la_1,\ad^m) \in \Co$
for all $m\geq 0, \la_1\prec \la.$ By Lemma \ref{escalar} again we
have
 $\e(\mu,\la_1)\e(\la_1,\mu) \in \Co$ for all $\mu\prec \ad^m,
\la_1\prec \la.$ Since by Lemma \ref{fusionrule} any $\mu$ with $
\col(\mu)= 0\mod \  n$ is a subsector of $\ad^m$ for some $m\geq 0$,
we conclude that $\e(\mu,\la_1)\e(\la_1,\mu) \in \Co$ for all $\mu,
\col(\mu)=0\mod \  n,\la_1\prec \la.$ By the definition of $S$
matrix we have
\[
|S_{\mu\la_1}|= d_{\la_1} |S_{\mu 1}|, \forall \mu, \col(\mu)=0 \mod
\ n\] Set $[a]=[b]=[1]$ in  Lemma \ref{orthophio} we have
\[
\l a_{\la_1}, a_{\la_1} \ra = \sum_{\mu,i,\col(\mu)= 0\mod \  n }
d_{\la_1}^2 \phi_1^{(\mu,i)}{\phi_1^{(\mu,i)}}^* = d_{\la_1}^2\] By
Lemma \ref{klx}  we have
\[ d_{\la_1}\geq \l a_{\la_1}, a_{\la_1}\ra \]
and we conclude that $d_{\la_1} =1,$ and in particular $v\la_1$ is
irreducible. The lemma now follows from Lemma \ref{fusionv}. \qed

The subfactors $a_\la(M)\subset M$ are type III analogue of
``orbifold subfactors" studied in \cite{EK} and \cite{Xo}.

\blem\label{no1} If $x\prec a_\la,$ $\la$ irreducible and $d_x=1$,
then $[\la]=[\omega^i],1\leq i\leq n$ and $[x]=[1]$. \elem

\prf If $[\la]\neq [\omega^i], \forall i,$ then by Lemma
\ref{fusionv} $\la\bar\la\succ \ad,$ and by Lemma \ref{irredad}  we
have $a_{\la}a_{\bar\la} \succ a_\ad$. Since $x\prec a_\la, d_x=1$,
by Lemma \ref{klx} we conclude that $d_{a_\ad}=d_\ad= 1.$ This is
impossible by Lemma \ref{irredad} and our assumption $k=n'n, n'\geq
3.$

\qed

Let $(n',n',...,n')$ be the unique fixed representation under the
action of $\Z_n.$ By Lemma \ref{klx}
\[[a_{(n',n'...,n')}] = \sum_{ 1\leq i\leq n} [b_i], [g^ib_1g^{-i}]=[b_{i+1}], 0\leq i\leq n-1 \]

\bdef\label{usub} Denote by $u:=(n'+1, n',n',...,n').$\ede Note that
by Lemma \ref{klx} $a_u$ is irreducible.

\blem\label{unot0} (1)
\[S_{u\ad}\neq 0;
\]
(2) Let $\Lambda=(n,0,...,0).$ Then $\l a_\Lambda, \tilde
a_{\bar\Lambda}\ra =0,$ and $S_{u\Lambda} \neq 0.$ \elem \prf Ad (1)
Since $n[a_u]=[a_v b_i],$ by Lemma \ref{orthophio}
\[ \frac{S_{u\ad}}{S_{1\ad}}=
\frac{S_{v\ad}}{nS_{1\ad}} \frac{S_{(n',...,n')\ad}}{S_{1\ad}}
\]

Direct computation using equation (\ref{kacp}) shows that
$\frac{S_{v\ad}}{S_{1\ad}}\neq 0.$ Note that by equation (\ref{ssy})
\[\frac{S_{(n',...,n')v}}{S_{1 v}} =0\] since $\col(v)= 1 $,
hence
\[
\frac{S_{(n',...,n')\ad}}{S_{1(n',...n')}}= -1
\]
and this implies that $S_{(n',...,n')\ad}\neq 0$ and (1) is
proved.\par Ad (2) Since $k=n'n\geq 3n,$ it follows that $\l\om^j
\Lambda, \bar\Lambda\ra =0, \forall 1\leq j\leq n.$ By Lemma
\ref{klx} $\l a_\Lambda, \tilde a_{\bar\Lambda}\ra =0.$ Since
$[a_va_{(n',n'...,n')}]= n [a_u],$ by Lemma \ref{orthophio} we have
\[
n \frac{S_{u\Lambda}}{S_{1\Lambda}} =
\frac{S_{v\Lambda}}{S_{1\Lambda}}\frac{S_{(n',...,n')\Lambda}}{S_{1\Lambda}}
\]
Hence to finish the proof we just have to check that
$S_{v\Lambda}\neq 0,S_{(n',...,n')\Lambda}\neq 0.$ Since
$Ch_{v'}(x_1,...,x_n)= \sum_{1\leq i\leq n} x_i,$ by equation
(\ref{kacp}) up to a nonzero constant $S_{v\Lambda}$ is equal to
\[
\exp(-2\pi i(2n-1)/(k+n)) + \sum_{ 0\leq j\leq n-2} \exp(-2\pi i
j/(k+n))
\]
This sum is equal to $0$ iff $n=k=2.$ Note that
$Ch_{\Lambda'}(x_1,...x_n)$ is a complete symmetric polynomial of
degree $n.$ $S_{v\Lambda}\neq 0$ now follows directly from equation
(\ref{kacp}) (cf. (2.7a) of \cite{GanW} for more general
statement).\qed

\qed

The main theorem of this section is:

\bthm\label{mn} The lattice of intermediate subfactors of
$a_u(M)\subset M$ is $M_{2n}.$ \ethm

The proof will be given in section \ref{proof}. Let us first show
that the subfactor in Theorem \ref{mn} contains $2n$ incomparable
intermediate subfactors. By fusion rule with $v$ in Lemma
\ref{fusionrule} we have
\[
[a_u]= [a_v b_i]=[b_i a_v], \forall 1\leq i\leq n. \] Therefore we
can assume that
\[ a_u= U_i a_vb_i U_i^* = V_ib_i a_v V_i^*, 1\leq i\leq n \]
where $U_i, V_i$ are unitaries. \bprop\label{2n} (1): As von Neumann
algebras \[U_i a_v(M) U_i^* = U_j a_v(M) U_j^*, V_i b_i(M) V_i^* =
V_j b_j(M) V_j^* \]
 iff $i=j$;\par (2) $U_ia_v(M)U_i^*$ is not an
intermediate subfactor in $V_jb_j(M)V_j^*\subset M$; \par
(3)$V_jb_j(M)V_j^*$ is not an intermediate subfactor in
$U_ia_v(M)U_i^*\subset M$. \eprop

\prf Ad(1): If $U_i a_v(M) U_i^* = U_j a_v(M) U_j^*,$ then $U_i
a_v(m) U_i^* = U_j a_v(\theta(m)) U_j^*, \forall m\in M$ where
$\theta$ is an automorphism of $M.$ By Frobenius reciprocity we have
$[\theta]\prec [a_va_{\bar v}]$. By Lemma \ref{no1} we conclude that
$[\theta]=[1]$ and hence \[U_i a_v(m) U_i^* = U_j a_v(U) a_v(m)
a_v(U)^* U_j^*, \forall m\in M \]
for some unitary $U\in M.$ Hence
\[ Ad_{U_i} a_v b_i  =Ad_{U_ja_v(U)} a_vb_i= Ad_{U_j} a_v b_j
\]
and we conclude that $[b_i]=[b_j]$, hence $i=j$. The second
statement in (1) is proved similarly. \par

Ad (2): If $U_ia_v(M)U_i^*$ is  an intermediate subfactor in
$V_jb_j(M)V_j^*\subset M$, then $Ad_{V_j} b_j = Ad_{U_i} a_v C$ for
some $C\in \End(M)$, and it follows that $[b_j\bar b_j]\succ
[a_v\bar a_v] \succ [a_\ad]$ Hence
\[
\l a_v b_j, a_v b_j \ra =\l b_j\bar b_j, a_v\bar a_v \ra \geq 2
\]
contradicting the irreducibility of $[a_u]=[a_vb_j]$. \par

Ad (3): If $V_jb_j(M)V_j^*$ is  an intermediate subfactor in
$U_ia_v(M)U_i^*\subset M,$ then there is $C'\in \End(M)$ such that
$[b_jC']=[a_v]$. Since $[a_v \bar a_v]=[1]+[a_\ad]$ and $a_\ad$ is
irreducible by Lemma \ref{irredad}, we must have $[b_j \bar
b_j]=[a_v\bar a_v]$ and therefore $d_{C'}=1$. By Frobenius
reciprocity $C'\prec [\bar b_j a_v],$ but $[\bar b_j a_v]$ is
irreducible since $a_u$ is irreducible, a contradiction. \qed

Here we give a quick proof of Th. \ref{mn} for $n=2$ and $k\neq 10,
28$ to illustrate some ideas behind the proof. Suppose that $M_1$ is
an intermediate subfactor of $ a_u(M)\subset M$. Since all factors
in this paper are isomorphic to hyperfinite type $III_1$ factor, we
can find $c_1, c_2\in \End(M)$ such that $a_u=c_1c_2$ and
$c_1(M)=M_1.$  Let $\rho=\rho_0 c_1$, and enumerate the basis of
$H_\rho$ by irreducible sectors $a$. Note that $a$ must be of the
form $\ro_0 c$ with $c$ irreducible, and so $d_a\geq d_{\ro_0}=\sqrt
2.$

Consider the fusion graph associated with the action of $v$ on
$H_\rho$: the vertices of this graph are irreducible sectors $a$,
and vertices $a$ and $b$ are connected by $\l va,b\ra$ edges. By
Lemma \ref{ephi} this graph is connected and has norm $2\cos
(\frac{\pi}{k+2})$, and hence it must be $A-D-E$ graph (cf. Chap. 1
of \cite{GHJ}). Since $k\neq 10,28$ it must be $A$ or $D$ graph. By
Lemma \ref{ephi} we have $\sum_a d_a^2 = \frac{1}{S_{11}^2} =
\frac{1}{\frac{1}{k+2} sin^2(\frac{\pi}{k+2})}.$  Since $d_a\geq
d_{\rho_0}= \sqrt 2,$ are the entries of Perron-Frobenius
eigenvector for the graph (Such eigenvector is unique up to a
positive scalar) , compare with the eigenvectors of A-D-E graphs
listed for example in Chap. 1 of \cite{GHJ}) we conclude that the
graph is $D$ graph and there is a sector $c$ with $d_c=1$ and
$c_1\prec a_\mu c$ for some $\mu\in \Delta$. We conclude that either
$[c_1]= [a_\mu c], $ or $[c_1]=[b_ic], 1\leq i\leq 2.$ In the former
case $[c_2]=[c^{-1} a_\la]$ or $[c_2]= [c^{-1} b_j], 1\leq j \leq
2$. But if $[c_2]=[c^{-1} a_\la]$ then $[a_u]=[a_\mu a_\la]$ is
irreducible, and by Lemma \ref{fusionv} $[a_\mu]=[a_u]$ or
$[a_\mu]=[1]$, which implies that $M_1$ is either $a_u(M)$ or $M.$
If $[c_2]= [c^{-1} b_j], 1\leq j \leq 2,$ then $[a_u]=[a_\mu b_j]$
and by computing the index and note that the colors of $u$ and $b_j$
are $ 1\mod \  2, 0\mod \  2$ respectively we have $a_\mu=a_v,$ and
we conclude that $M_1$ must be one of the intermediate subfactors
given in Prop. \ref{2n}. The case of $[c_1]=[b_ic], 1\leq i\leq 2$
is treated similarly. By Prop. \ref{2n} we have proved Th.\ref{mn}
for $n=2, k\neq 10,28.$ The same idea as presented above can be used
to give a complete list of all intermediate subfactors of
Goodman-Harpe-Jones subfactors. We hope to discuss this and related
problems elsewhere.

\section{Centrality of a class of intertwinners and its
consequences}\label{ucentral} We preserve the setup of section
\ref{typea}.

Assume that $\ro\bar\ro\in \Delta_\A.$ We will investigate a class
of inductions which are motivated by finding a proof of Th.
\ref{mn}.

In this section we assume that $[a_v]=[h\tilde a_v], [h^n]=[1],$
$a_{\ad}$ is irreducible, and if $\mu\prec \ad^2, 1\prec a_\mu,$
then $[\mu]=[1].$

Choose a unitary $T\in \Hom (a_v,h\tilde a_v).$ Such $T$ is unique
up to scalar since $a_v$ is irreducible. By Lemma  \ref{acommut} we
have $[h \tilde a_v]=[\tilde a_vh].$ Choose a unitary $T_1\in
\Hom(\tilde a_v h, h\tilde a_v).$ Note that $T_1$ is unique up to
scalar since $h\tilde a_v$ is irreducible.

\bdef\label{Ts} Denote by $U_n:=Ta_v(T) a_v^2(T)...a_v^{n-1}(T)\in
\Hom (a_{v}^n, (h\tilde a_v)^n). $

Denote by $T_i:=T_1\tilde a_v (T_1)...\tilde a_v^{i-1}(T_1) \in \Hom
(\tilde a_v^i h, h\tilde a_v^i), 1\leq i\leq n-1.$

Choose $T'\in \Hom (h^n,1)$ ($T'$ is unique up to scalar). \ede

\bdef\label{uw} Set $w=v^n$ and define $u_w:= T' h^{n-1}(T_{n-1})
h^{n-2}(T_{n-2})... h(T_1)U_n \in \Hom (a_v^n, \tilde a_v^n).$\ede

For example when $n=3,$ $u_w= T' h^2(T_1) h^2(\tilde a_v(T_1))h(T_1)
T a_v(T) a_v^2(T).$ The reader is encouraged to give a diagrammatic
representation of $u_w$ as in \cite{Xb}.

\blem\label{c} Suppose that $x,y$ are sectors such that
$$[x]=\sum_{1\leq i\leq m} [x_i], [y]=\sum_{1\leq i\leq m} [y_i],
d_{x_i}< d_{x_j}, d_{y_i}< d_{y_j}$$
if $i<j,$ and $x_i,y_i $ are
irreducible.  Let $T_{x,i}\in \Hom( x_i,x), T_{y,i}\in \Hom(y_i,y),
i=1,...,m$ be isometries.

If $U\in \Hom(x,y)$ is unitary then $UT_{x,i} T_{x,i}^* U^*=
T_{y,i}T_{y,i}^*,i=1,...,m.$ \elem

\prf By assumption $\Hom(x,x), \Hom(y,y)$ are finite dimensional
abelian algebras, and so for each $1\leq i\leq m$ we have $UT_{x,i}
T_{x,i}^* U*= T_{y,j}T_{y,j}^*$ for some $j.$

By equation (\ref{phiprop}) we have
\[
d_y \phi_y(UT_{x,i} T_{x,i}^* U^*)= d_x \phi_x(T_{x,i} T_{x,i}^* ) =
d_{x_i}
\]
Hence $d_{x_i}= d_{y_j}.$ By assumption it follows that $i=j, 1\leq
i\leq m.$ \qed

\blem\label{c0} Let $U\in \Hom (a_v^2 h^j, h^i \tilde a_v^2),
i,j\geq 0$ be a unitary. Then $h^i (\bar\ro(\e(v,v))) U= U
\bar\ro(\e(v,v)).$ \elem

\prf Since $a_\ad$ is irreducible, we have $\l a_va_v, a_va_v\ra =\l
a_v\bar a_v,a_v\bar a_v\ra =2.$  We note that
$[a_va_v]=[a_{(2,0,...,0)}]+[a_{(0,1,0,...,0)}]$ and
$\frac{d_{a_{(2,0,...,0)}}}{d_{a_{(2,0,...,0)}}}=
\frac{sin(\frac{(n+1)\pi}{k+n})} {sin(\frac{(n-1)\pi}{k+n})} >1 $
and so the assumption of Lemma \ref{c} is verified. Denote by
$P_1,P_2\in \Hom (v^2,v^2)$ the  two different minimal projections
corresponding to $(2,0,...,0), (0,1,...,0)$ respectively. Note that
$\bar\ro(P_l), h^i(\bar\ro(P_l)), l=1,2$ are minimal projections in
$\Hom (a_v^2 h^j,a_v^2 h^j), \Hom(h^i \tilde a_v^2,h^i \tilde
a_v^2)$ respectively and by Lemma \ref{c} we have
$U^*h^i(\bar\ro(P_l)) U= \bar\ro(P_l), l=1,2.$

Assume that $\e(v,v)= z_1 P_1 +z_2 P_2$ where $z_1,z_2\in \Co$ (cf.
Lemma 3.1.1 \cite{Xj} for explicit formulas for $z_1,z_2$). Then
$h^i(\bar\ro(\e(v,v)))= z_1 h^i(\bar\ro(P_1)) + z_2
h^i(\bar\ro(P_2))$ and the lemma follows. \qed

\blem\label{c1} $\tilde a_v^i (\bar\rho(\e(v,v))) u_w = u_w a_v^i
(\bar\ro(\e(v,v)), 0\leq i\leq n-2.$ \elem

\prf By Def. \ref{uw} we can write $u_w= V_1'V_2'V_3'$  where
\[V_3'= a_v^{i+2}(V_3), V_3= h^{n-i-3}(T_{n-i-3})...h^2(T_2)
h(T_1)\in \Hom (a_v^{n-i-2}, h^{n-i-2} \tilde a_v^{n-i-2})\]

$V_2'= a_v^i (V_2), V_2= h^{n-i-1}(T_2)...h^2(T_2) h(T_1) T
a_v(T)\in \Hom (a_v^2 h^{n-i-2}, h^{n-i} \tilde a_v^2) $ and $V_1'=
T' h^{n-1}(T_i)...h^i(T_i) h^{i-1}(T_{i-1}) h^{i-2}(T_{i-2})...
h(T_1) T a_v(T)... a_v^{i-1}(T) \in \Hom (a_v^i h^{n-i}, \tilde
a_v^i).$

Although the complicated but explicit formulas of $V_1', V_2, V_3$
are given above,  we only use their intertwining properties in what
follows.

Hence \begin{align*} \tilde a_v^i(\bar\ro(\e(v,v))) u_w &=\tilde
a_v^i(\bar\ro(\e(v,v))) V_1' a_v^i(V_2) a_v^{i+2}(V_3) \\
&= V_1' a_v^i (h^{n-i}(\bar\ro(\e(v,v))) V_2)a_v^{i+2}(V_3)  \\
&= V_1' a_v^i(V_2 \bar\ro(\e(v,v)))a_v^{i+2}(V_3)\\
&=V_1' a_v^i(V_2) a_v^i(\bar\ro(\e(v,v))a_v^{2}(V_3))\\
&=V_1' a_v^i(V_2) a_v^{i+2}(V_3) a_v^i(\bar\ro(\e(v,v)))\\
&= u_wa_v^i(\bar\ro(\e(v,v)))
\end{align*}
where in the third $=$ we have used Lemma \ref{c0}. \qed

\blem\label{c2} $\tilde a_v^{n-1}(\bar\ro(\e(v,v))) u_w  a_w(u_w)
=u_w a_w(u_w) a_v^{n-1}(\bar\ro(\e(v,v))).$ \elem

\prf By Def. \ref{uw} we can write $u_w a_w(u_w)= W_1'W_2'W_3'$
where $W_3' =a_v^{n+1}(W_3), W_3 = h^{n-2}(T_{n-2})...h(T_2)h(T_1) T
a_v(T)... a_v^{n-2}(T) \in \Hom(a_v^{n-1}, h^{n-1}\tilde a_v^{n-1})$

$W_2'= a_v^{n-1}(W_2), W_2 =TT' h^{n-1}(T_1)...h(T_1) a_v(T)\in
\Hom(a_v^2 h^{n-1}, h\tilde a_v^2)$

and

$ W_1'= T' h^{n-1}(T_n) h^{n-2}(T_{n-2}) ... h(T_1) T
a_v(T)...a_v^{n-2}(T)\in \Hom (a_v^{n-1}h, \tilde a_v^{n-1})$

As in the proof of Lemma \ref{c1},  even though  explicit formulas
of $W_2, W_3, W_1'$ are given as above, what we need in the
following is their intertwining properties.

Hence
\begin{align*}
\tilde a_v^{n-1}(\bar\ro(\e(v,v))) u_w \tilde a_w(u_w) &=
\tilde a_v^{n-1}(\bar\ro(\e(v,v))) W_1' a_v^{n-1}(W_2) a_v^{n+1}(W_3) \\
&= W_1' a_v^{n-1} (h(\bar\ro(\e(v,v))) W_2)a_v^{n+1}(W_3) \\
&= W_1' a_v^{n-1}(W_2)a_v^{n-1} (\bar\ro(\e(v,v)) a_v^2(W_3))\\
&=W_1' a_v^{n-1}(W_2) a_v^{n+1}(W_3)a_v^{n-1} (\bar\ro(\e(v,v)) )\\
&=u_w  a_w(u_w) a_v^{n-1}(\bar\ro(\e(v,v)))
\end{align*}

Where in the third $=$ we have used Lemma \ref{c0}. \qed

\bdef\label{uwm} For each integer $m\geq 1,$
$u_{w^m}:=u_wa_w(u_w)... a_w^{m-1} (u_w) \in \Hom(a_{w^m}, \tilde
a_{w^m}).$ \ede

\bthm\label{centu} Let $m\geq 1$ be any integer and $R\in \Hom( w^m,
w^m).$ Then
$$
\bar\ro(R) u_{w^m}=u_{w^m}\bar\ro(R).
$$\ethm
\prf By Lemma \ref{denseb} it is sufficient to prove the theorem for
$R= v^{m'}(\bar\ro(\e(v,v))), 1\leq m'\leq m-1.$ When $nn_1< m'<n(
n_1+1),n_1\in \Z$ we can write
$$
u_wa_w(u_w)... a_w^{m-1} (u_w) =U_1' a_w^{n_1}(u_w) U_2'$$ where
$U_1'\in \Hom(a_w^{n_1}, a_w^{n_1}), U_2'\in a_w^{n_1+1}(M)$ and the
theorem follows from Lemma \ref{c1}. Similarly when $m'=nn_1, n_1\in
\Z$ we can write
$$
u_wa_w(u_w)... a_w^{m-1} (u_w) =U_1'' a_w^{n_1-1}(u_w a_w(u_w))
U_2''$$ with $U_1''\in \Hom(a_w^{n_1-1}, a_w^{n_1-1}), U_2''\in
a_w^{n_1+2}(M)$ and the theorem follows from Lemma \ref{c2}. \qed

\blem\label{indep} Suppose that $\mu\prec w^m $ are irreducible and
Let $t_{\mu,i}\in \Hom(\mu, w^m), m\geq 1$ be a set of isometries
such that $\sum_{\mu,i} t_{\mu,i}t_{\mu,i}^* =1.$ Then \par (1)For
each fixed $\mu$, $\bar\ro(t_{\mu,i})^* u_{w^m}
\bar\ro(t_{\mu,i})\in \Hom(a_\mu,\tilde a_\mu)$  is independent of
choices of $t_{\mu,i};$
\par (2)$\bar\ro(t_{\mu,i})^* u_{w^m} \bar\ro(t_{\mu,i})\in
\Hom(a_\mu,\tilde a_\mu)$ is unitary. \elem

\prf (1) follows immediately from Th. \ref{centu}. To prove (2),
note that for each fixed $\mu,i$
\[
1=\sum_{\la,j}\bar\ro(t_{\mu,i})^* u_{w^m} \bar\ro(t_{\la,j})
\bar\ro(t_{\la,j})^* u_{w^m}^* \bar\ro(t_{\mu,i})^*
=\bar\ro(t_{\mu,i})^* u_{w^m} \bar\ro(t_{\mu,i})\bar\ro(t_{\mu,i})^*
u_{w^m}^* \bar\ro(t_{\mu,i})\] where in the second $=$ we have used
Th. \ref{centu}. Similarly
\[1=\bar\ro(t_{\mu,i})^* u_{w^m}^* \bar\ro(t_{\mu,i})\bar\ro(t_{\mu,i})^*
u_{w^m} \bar\ro(t_{\mu,i})\] and the Prop. is proved.\qed

The unitary in (2) of Prop. \ref{indep} will be denoted by $u_\mu$
(it may depend on $m$) in the following.

 \bdef\label{psiw} Let $\mu\in \Delta_\A$ and $b\in H_\ro$ be
 irreducible. Define
$$\psi_b^{(w)}:= S_{11} d_b d_w\phi_w (\e(b\bar\ro,w) b(u_w)
\e(w,b\bar\ro)), b\in H_\ro.$$ \ede

\blem\label{sumpsi} Let $m\geq 1$ $t_{\mu,i}$ be as in Prop.
\ref{indep}. Then
$$
|\sum_b d_b^2(\frac{\psi^{(w)}_b}{d_b S_{11}})^m| =
\frac{1}{S_{11}^2} \l w^m,1\ra , \forall m\geq 1.$$
\elem

\prf
\begin{align*}
 (\frac{\psi^{(w)}_b}{d_b S_{11}})^m &= d_w^m \phi_w^m
(\e(b\bar\ro,w^m) b(u_{w^m}) \e(w^m,b\bar\ro))\\
&= \sum_{\mu,i} d_\mu \phi_\mu(t_{\mu,i}^* \e(b,w^m)b(u_{w^m})
\e(w^m,b\bar\ro)t_{\mu,i}) \\
&=\sum_{\mu,i} d_\mu \phi_\mu(\e(b\bar\ro,\mu)b(\bar\ro(t_{\mu,i})^*
u_{w^m}\bar\ro(t_{\mu,i})) \e(\mu,b\bar\ro))\\
&= \sum_\mu \l \mu, w^m\ra d_\mu \phi_\mu(\e(b\bar\ro,\mu)u_\mu
\e(\mu,b\bar\ro))
\end{align*}
where we have used definition of minimal left inverse in the first
$=$, equation (\ref{phiprop}) in the second $=$, Prop. \ref{bfe} in
the third $=$, and Lemma \ref{indep} in the last $=$.

It follows that

\begin{align*}
 \sum_b d_b^2(\frac{\psi^{(w)}_b}{d_b S_{11}})^m &= \sum_{b,\mu}
\l \mu, w^m\ra d_b^2 d_\mu \phi_\mu(\e(b\bar\ro,\mu)u_\mu \e(\mu,b\bar\ro)) \\
&= \sum_{\mu}\l \mu,w^m\ra d_\mu \sum_{b} d_b
d_b\phi_\mu(\e(b\bar\ro,\mu)b(u_\mu)\e(\mu,b\bar\ro)) \\
&=\sum_b d_b^2 \phi_1(u_1) \l 1,w^m\ra
\end{align*}

where we have used Lemma \ref{epsi} in the third $=$. Since $u_1 \in
\Hom (1,1)$ is unitary by Prop. \ref{indep}, $|\phi_1(u_1)|=1$ and
we have proved that
$$
|\sum_b d_b^2 (\frac{\psi^{(w)}_b}{d_b S_{11}})^m| =
\frac{1}{S_{11}^2} \l w^m,1\ra.
$$

\bprop\label{extrempsiw} There is a sector $c\in H_\ro$ such that
$|\frac{\psi_c^{(w)}}{S_{11}}|= d_c d_w.$ \eprop

\prf By Lemma \ref{sumpsi} we have
$$
|\sum_b d_b^2(\frac{\psi^{(w)}_b}{d_b S_{11}})^m| =
\frac{1}{S_{11}^2} \l w^m,1\ra, \forall m\geq 1.
$$
By repeated using Verlinde formula we have
$$
\l w^m,1\ra = \sum_\mu \frac{1}{S_{1\mu}^2}
(\frac{S_{v\mu}}{S_{1\mu}})^{nm}
$$
By Lemma \ref{scalarv},  when $m$ goes infinity, the leading order
of $|\sum_b d_b^2(\frac{\psi^{(w)}_b}{d_b S_{11}})^m| $must be $n
d_w^m$. Note by Lemma \ref{epsi} $|\frac{\psi^{(w)}_b}{d_b
S_{11}}|\leq d_w.$ It follows that There is a sector $c\in H_\ro$
such that $|\frac{\psi_c^{(w)}}{S_{11}}|= d_c d_w.$ \qed

Choose $m=1$  and let $t_{\mu,i}$ be isometries as in Lemma
\ref{indep}.

\bdef\label{psimu} Assume that $\mu\in \Delta_\A$ and $[b]\in H_\ro$
is irreducible. Define
\[
\frac{\psi_b^{(\mu)}}{S_{11}}:=d_b d_\mu \phi_\mu ((\e(b\bar\ro,\mu)
b(\bar\ro(t_{\mu,i})^* u_w \bar\ro(t_{\mu,i}))\e(\mu,b\bar\ro))).
\]
\ede Note that by Lemma \ref{indep} $\psi_b^{(\mu)}$ is independent
of the choice of $i$.

\bcor\label{smallc} Assume that $[a_v]=[h\tilde a_v], [h^n]=[1],$
$a_{\ad}$ is irreducible, and if $\mu\prec \ad^2, 1\prec a_\mu,$
then $[\mu]=[1].$ Then there is $[c]\in H_\rho$ such that
$|\frac{\psi_c^{(\la)}}{S_11}|=d_c d_\la, \forall \la, \col(\la)= 0
\mod \ n$ and $[c\bar c]=\sum_{1\leq i_2\leq \frac{n}{i_1}}
[\om^{i_2i_1}] $ where $i_1$ is a divisor of $n$. \ecor

\prf Choose $m=1$  and let $t_{\mu,i}$ be isometries as in Lemma
\ref{indep}. By equation (\ref{phiprop}) we have
\[
\frac{\psi_c^{(w)}}{S_{11}}= \sum_{\mu}\l \mu,w\ra
\frac{\psi_c^{(\mu)}}{S_{11}}
\]
By Lemma \ref{epsi} we have
$$
|\frac{\psi_c^{(\mu)}}{S_{11}}|\leq d_cd_\mu
$$
By Prop. \ref{extrempsiw} we conclude that
$$
|\frac{\psi_c^{(\mu)}}{S_{11}}|= d_cd_\mu ,\forall \mu\prec w$$ In
particular $|\frac{\psi_c^{(\ad)}}{S_{11}}|= d_cd_\ad .$ By Lemma
\ref{epsi} we know that $\sum_b {\psi_b^{(\ad)}}^* [b]$ is a nonzero
eigenvector of the action of $[\la]$ on $H_\ro$. Since $\l
a_\ad,\tilde a_\ad\ra =1,$ by Prop. \ref{spec} we must have
$\psi_b^{(\ad)} = z \phi_b^{(\ad)}, $ for some constant $z$
independent of $b.$  Since $[\bar\ad]=[\ad]$, $\sum_b
\phi_b^{(\ad)}b $ is also an eigenvector of the action of $[\la]$
with eigenvalue $\frac{S_{\la\ad}}{S_{1\ad}},$ it follows that
${\phi_b^{(\ad)}}^* = z' \phi_b^{(\ad)}, $ for some constant
$|z'|=1$ independent of $b.$ Hence
\[
\sum_b{\psi_b^{(\ad)}}^2 = \sum_b z^2 \bar z'\phi_b^{(\ad)}
{\phi_b^{(\ad)}}^* = z^2 \bar z'
\]
By (3) of Lemma \ref{epsi} and our assumption we conclude that
$|z|=1$, and so by Lemma \ref{ephi} we have

\[
d_c^2 =|\frac{\psi_c^{(\ad)}}{S_{1\ad}}|^2 =
|\frac{\phi_c^{(\ad)}}{S_{1\ad}}|^2=\sum_\mu \l c\bar c, \mu\ra
\frac{S_{\mu\ad}}{S_{1\ad}}
\]
Since $\frac{S_{\mu\ad}}{S_{1\ad}}\leq d_\mu$, we must have
$\frac{S_{\mu\ad}}{S_{1\ad}}= d_\mu, \forall \mu\prec  c\bar c.$

By Lemma \ref{scalarad} we conclude that if $\mu\prec c\bar c $,
then $\mu =\om^i$ for some $1\leq i\leq n.$  Let $1\leq i_1\leq n$
be the smallest positive integer such that $[\om^{i_1}c]=[c].$ Then
it is clear that $[c\bar c]=\sum_{1\leq i_2\leq \frac{n}{i_1}}
[\om^{i_2i_1}] $ where $i_1$ is a divisor of $n$. \qed
\section{Proof of Th. \ref{mn}}\label{proof}
In this section we preserve the setting of section \ref{orb}.  Let
$c_1, c_2 \in \End(M)$ such that $a_u =c_1 c_2, c_1(M)=M_1, M_1\neq
a_u(M), M$. By Prop. \ref{2n} to prove Th. \ref{mn} it is enough to
show that $M_1$ is one of the intermediate subfactors in Prop.
\ref{2n}.
\subsection{Local consideration}\label{localconsideration}
Suppose $c$ is a sector such that $c\bar c \prec a_\mu^{\ro_o}$
where $\mu\in \Delta_\A$  is a direct sum of irreducible sectors
with colors divisible by $n.$ Recall from section \ref{orb} that if
$\la=0 \mod \ n,$ then $[a^{\ro_o}_\la]=[\tilde a^{\ro_o}_\la],$ and
we can apply induction of $a^{\ro_o}_\la$ with respect to $c.$ The
following Lemma is proved by  a translation of the proof of (3) of
Lemma 3.3 in \cite{Xp} into our setting: \blem\label{applyxp} If
$\la=0 \mod \ n,$ then $[a_{a_\la^{\ro_o}}^{c}] =[a_\la^{\ro_o c}].$
\elem

By Prop. \ref{loclr} we have $c_1=c_1'c_1''.$ Let $c_2'=c_1''c_2$ so
that $a_u=c_1'c_2'.$ Consider induction with respect to $\ro_o
c_1'.$

We have

\blem\label{c'=1} $[c_1'\bar c_1']=[1].$\elem

\prf Apply Lemma \ref{ephi} to $a=\ro_0 c_1', b=\ro_0 \bar c_2'$ we
have
\[\sum_{i} \frac{\phi_a^{(\la,i)}{\phi_b^{(\la,i)}}^*}{S_{1\la}^2}=
\sum_\nu \l \ro_0 c_1'c_2'\bar\ro_0, \nu\ra
\frac{S_{\nu\la}}{S_{1\la}} =\sum_\nu \l u\ro_0 \bar\ro_0, \nu\ra
\frac{S_{\nu\la}}{S_{1\la}} =\sum_{1\leq i\leq n} \exp(\frac{2\pi i
col(\la)}{n}) \frac{S_{u\la}}{S_{1\la}}
\]
Choose $\la=\ad$ and use Lemma \ref{unot0} we have
\[\sum_{i} \frac{\phi_a^{(\la,i)}{\phi_b^{(\la,i)}}^*}{S_{1\la}^2}\neq
0\] Hence by Lemma. \ref{znot0} we obtain $\l a_\ad^{\ro_0c_1'},
\tilde a_\ad^{\ro_o c_1'}\ra\geq 1.$ For any $\mu\in \Delta_\A,$
since $\ro_0 c_1' \bar c_1' \bar\ro_o\prec \ro_o a_{u\bar u}
\bar\ro_o $ and each irreducible sector of $[\ro_o a_{u\bar u}
\bar\ro_o ]= [\ro_0 \bar \ro_o u\bar u]$ has color divisible by $n,$
it follows that if $\col (\mu) \neq 0 \mod \ n,$ then $\l \mu, \ro_0
c_1' \bar c_1' \bar\ro_o \ra =0.$ On the other hand if  $\col (\mu)
=0 \mod \ n,$ by Lemma \ref{applyxp} and Prop. \ref{xua} we have
\[ \l a_\mu^{\ro_o c_1'}, 1\ra= \l a_\mu^{\ro_o}, c_1'\bar c_1' \ra
= \l \mu, \ro_0 c_1' \bar c_1' \bar\ro_o \ra
\]
By (1) of Prop. \ref{loclr} it follows that $\ro_o c_1'$ is local.

 By Lemma \ref{irredad}
we have $[a_\ad^{\ro_0c_1'}] =[\tilde a_\ad^{\ro_0c_1'}]$, and by
Lemma \ref{loc2} and Lemma \ref{scalarad} we conclude that $[\ro_o
c_1'\bar c_1' \bar\ro_o] =\sum_j [\omega^j]$ where the sum is over a
finite set of positive integers. Since $\ro_0 c_1'$ is irreducible
and $[\ro_0\bar\ro_0 ]=\sum_{1\leq j\leq n} [\omega^j]$ we conclude
that $[\ro_0 c_1'\bar c_1' \bar\ro_0] =\sum_{1\leq j\leq n}
[\omega^j]$ Hence $d_{c_1'}=1$ and $[c_1'\bar c_1'] =[1].$ \qed

By Prop. \ref{loclr} we have proved

\bcor\label{cor'} If $\la\in \Delta_\A$ is irreducible,  then $\l 1,
a^{\ro_o c_1}_\la\ra \geq 1 $ iff $ \la = \omega^i , 1\leq i \leq
n.$ \ecor
\subsection{Verifying assumptions of
Cor. \ref{smallc}}\label{verify} Set $\ro= \ro_0 c_1$ and all
inductions in the rest of this section are with respect to $\rho$.

\blem\label{irredad2} $a_\la$ is irreducible for all irreducible
descendants of $v^2\bar v^2, v\bar v^3.$ \elem

\prf By Lemma \ref{fusionrule} and Prop. \ref{xua} we have for
$n\geq 3$
\[
[a_{v\bar v}a_{v\bar v}]= 2[1]+ 4[a_\ad] + [a_{(2,0,...,0,2)}] +
[a_{(0,1,0,...,1,0)}]+ [a_{(0,1,0,...,0,2)}] + [a_{(2,0,...,0,1,0)}]
\]
Note that by Cor. \ref{cor'}  we have
\[\l a_\la, a_\mu\ra  = \l 1, a_{\bar\la \mu}\ra\geq 2\]
iff $[\om^j(\la)]=[\mu]$ for some $1\leq j\leq n-1.$ It is easy to
check with the explicit formulas above that  $a_\la$ is irreducible
for all irreducible descendants of $v^2\bar v^2.$ $n=2$ case is
simpler, and similarly one can check directly that  $a_\la$ is
irreducible for all irreducible descendants of $ v\bar v \bar v^3.$
\qed

\blem\label{nobar} For all $\la$ with $\col(\la)=0,[a_\la]= [\tilde
a_\la].$ \elem \prf By (2) of Prop. \ref{spec} and Th. 2.1 of
\cite{Gan} all $Z_{\la\mu}$ with $Z_{1,\la}\neq 0$ iff $\la=\om^i
,1\leq i\leq n$ are classified. Using Cor. \ref{cor'} , it follows
by inspection of Th. 2.1 of \cite{Gan} that for all $\la$ with
$\col(\la)=0,$ $Z_{\la\la}= \l a_\la, \tilde a_\la\ra \neq 0$ or
$Z_{\la\la}= \l a_{\bar \la}, \tilde a_\la\ra \neq 0, \forall \la.$
In the latter case by Prop. \ref{spec} we conclude that $\la$
appears in $\Exp$ iff $\l a_\la, a_{\bar\la}\ra \neq 0.$ Choose
$\la=(n,0,...,0)=\Lambda$ as in Lemma \ref{unot0}. It follows from
Lemma \ref{unot0} and Cor. \ref{znot0} that $\Lambda \in \Exp,$ but
$\l a_\Lambda, \bar a_\Lambda\ra =0,$ contradiction. Hence $\l
a_\la, \tilde a_\la\ra \neq 0, \forall \la, \col(\la)=0\mod \  n,$
and by Lemma \ref{klx} we conclude that for all $\la$ with
$\col(\la)=0,$ $[a_\la]= [\tilde a_\la].$ \qed \blem\label{xy}
Suppose that $x_i\prec a_{\la_i}\tilde a_{\mu_i}, i=1,2$ and
$x_1x_2$ is a direct sum of $a_\nu$ with $a_\nu$ irreducible. Then
$[x_1x_2]=[x_2x_1].$ \elem

\prf By assumption it is enough to check that
\[\l x_1x_2, a_\nu\ra =\l x_2x_1, a_\nu\ra \]
By Lemma  \ref{acommut} we have $[a_\nu \bar x_2]=[\bar x_2 a_\nu],$
together with Frobenius reciprocity we obtain
\[\l x_1x_2, a_\nu\ra= \l x_1, a_\nu\bar x_2\ra= \l x_1, \bar
x_2a_\nu\ra\ra =\l x_2x_1, a_\nu\ra \] \qed

 \bprop\label{lemh}
There exists $h\in \End(M)$ such that $[\tilde a_v]=[h a_v],
[h^n]=[1].$ \eprop

\prf First suppose that there is no automorphism $h$ such that
$[\tilde a_v] =[h a_v]$ or $[\tilde a_{\bar v}] =[h a_v].$  By Lemma
\ref{nobar} $[a_va_{\bar v}]=[\tilde a_v\tilde a_{\bar
v}]=[1]+[a_\ad].$ By Lemma \ref{irredad} $a_\ad$ is irreducible, it
follows that there are sectors $x_i,y_i$ with $d_{x_i}>1, d_{y_i}>1$
such that
\[
[a_v\tilde a_v]= [x_1]+[x_2], [a_{\bar v}\tilde a_v]= [y_1]+[y_2].
\]
We compute
\[ [a_v\tilde a_v\tilde a_{\bar v}]= [x_1 \tilde a_{\bar v}] + [x_2 \tilde a_{\bar v}]
=[a_v a_v a_{\bar v}]= 2[a_v]+
[a_{(2,0,...,0,1)}]+[a_{(0,1,0,...,0,1)}]
\]
where we have used Lemma. \ref{nobar} in the second $=.$ By
assumption $d_{x_i}> 1, i=1,2$ we have $x_i \tilde a_{\bar v}\succ
a_v,$ but $[x_i \tilde a_{\bar v}]\neq  [a_v], i=1,2.$ Hence we can
assume that
\[
[x_1 \tilde a_{\bar v}]= [a_v]+[a_{(2,0,...,0,1)}], [x_2 \tilde
a_{\bar v}]=[a_v]+[a_{(0,1,0,...,0,1)}]\] Hence
\[ \l a_{\bar v} x_i, a_{\bar v} x_i\ra = \l x_i a_{\bar v}, x_i a_{\bar v} \ra
=\l x_i , x_i a_{\bar v v} \ra =\l x_i , x_i \tilde a_{\bar v v} \ra
=2
\]
where we have used Lemma \ref{acommut} in the first $=$ and Lemma.
\ref{nobar} in the third $=.$ We can assume that
\[
[a_{\bar v} x_i]= [\tilde a_v]+[u_i], i=1,2
\]
where $u_i, i=1,2$ is irreducible and we may have $[u_1]=[u_2].$
Note that $[a_{\bar v}x_1]+ [a_{\bar v}x_2]=[a_{ v}y_1]+ [a_{
v}y_2]= [a_{\bar v}a_v \tilde a_v].$

The same argument applies to $y_i, i=1,2$ and we may choose $y_i$
such that
\[
[a_{\bar v} x_i]= [a_v y_i], i=1,2
\]
Consider now
\begin{align*}
[a_{v\bar v}^2] & = [x_1 \bar x_1] +[ x_2 \bar x_2]+[x_1 \bar
x_2]+[x_2 \bar x_1] \\
&=2[1]+ 4[a_\ad] + [a_{(2,0,...,0,2)}] + [a_{(0,1,0,...,1,0)}]+
[a_{(0,1,0,...,0,2)}] + [a_{(2,0,...,0,1,0)}]
\end{align*}
Note that $x_i\bar x_i\succ a_{v\bar v},$ and $[x_i\bar x_j]=[\bar
x_jx_i]$ by Lemma \ref{irredad2} and Lemma \ref{xy}. Hence
\[
\l x_2\bar x_1, x_2\bar x_1\ra= \l x_2\bar x_2, x_1\bar x_1\ra \geq
2
\]
By computing the index of sectors we conclude that
\begin{align*}
[x_1 \bar x_1]=[a_{v\bar v}] + [a_{(2,0,...,0,2)}], &[x_1 \bar
x_2]=[a_{\ad}] + [a_{(0,1,...,0,2)}] \\
[x_2 \bar x_2]=[a_{v\bar v}] + [a_{(0,1,0,...,1,0)}], &[x_2 \bar
x_1]=[a_{\ad}] + [a_{(2,0,...,0,1,0)}]
\end{align*}
Similarly we obtain
\begin{align*}
[y_1 \bar y_1]=[a_{v\bar v}] + [a_{(2,0,...,0,2)}], &[y_1 \bar
y_2]=[a_{\ad}] + [a_{(0,1,...,0,2)}] \\
[y_2 \bar y_2]=[a_{v\bar v}] + [a_{(0,1,0,...,1,0)}], &[y_2 \bar
y_1]=[a_{\ad}] + [a_{(2,0,...,0,1,0)}]
\end{align*}

Next compute
\[
[a_{\bar v^2} a_{v\bar v}]= [a_{\bar v} \tilde a_v a_{\bar v} \tilde
a_{\bar v}] = [y_1\bar x_1]+ [y_1\bar x_2]+[y_2\bar x_1]+[y_2\bar
x_2].
\]
Note that
\[
\l y_2\bar y_1, x_2\bar x_1\ra =\l \bar x_2 y_2,\bar x_1 y_1\ra =2
\]
\[
2=\l a_{\bar v} x_i, a_v y_i\ra = \l a_{\bar v}^2, y_i \bar x_i\ra
\]
and
\[
\l y_i\bar x_i, y_i\bar x_i\ra =\l y_i \bar y_i, x_i\bar x_i\ra =3
\]
\[
\l y_1\bar x_2, y_1 \bar x_2\ra = \l y_1\bar y_2, x_2 \bar x_2\ra=2
\]
where  we have also used Lemma \ref{xy}. From these equations we
conclude that
\[
[y_1\bar x_1] = [a_{\bar v}^2] + [a_{(1,0,...,0,3)}]
\]
or
\[
[y_1\bar x_1] = [a_{\bar v}^2] + [a_{(1,0,...0,1,0,0)}]
\]
From $[a_{\bar v}x_1]=[a_vy_1]$ we obtain
\[
[a_{\bar v} x_1\bar x_1]= [a_vy_1\bar x_1]
\]
Using the formulas for $x_1\bar x_1, y_1\bar x_1$ we obtain
\[
[a_{\bar v} a_{(2,0,...,0,2)}] =[ a_v a_{(1,0,...,0,1,0,0)}]
\]
or
\[
[a_{\bar v} a_{(2,0,...,0,2)}] =[ a_v a_{(1,0,...,0,3)}]
\]
Both identities are incompatible with Lemma \ref{fusionrule} and
Lemma \ref{nobar}.\par Therefore there is an automorphism $h$ such
that $[\tilde a_v] =[h a_v]$ or $[\tilde a_{\bar v}] =[h a_v].$
Hence $h^n\prec [\tilde a_{\bar v^n}a_{v^n}]= [a_{\bar v^n}a_{v^n}]$
or $h^n\prec [\tilde a_{v^n}a_{v^n}]=[a_{ v^n}a_{v^n}]$ by Lemma
\ref{nobar}. Assume that $h^n\prec a_\mu$ for some $\mu,
\col(\mu)=0\mod \  n.$  Since $\ro=\ro_0c_1,$  by Lemma
\ref{applyxp} there is a sector $x$ of $a_\mu^{\ro_0}$ such that
$[a_x^{c_1}]=[h^n].$  Since $d_x=1,$ by Lemma \ref{no1} we conclude
that $[x]=[1]$ and  $[h^n]=[1].$ \par If $[\tilde a_{\bar v}] =[h
a_v],$ use $[h^n]=[1]$ we have $[a_{v^n}]=[a_{\bar v^n}]$ Hence
$\om^j (n,0,...,0) \prec \bar v^n$ for some $1\leq j\leq n$ which is
incompatible with fusion rules in Lemma \ref{fusionrule} since
$k=n'n\geq 3n.$ \qed
\subsection{Properties of sectors related to $a_u$}\label{relatedtou}

\blem\label{holo} If $\e(\om^l,\la)\e(\la,\om^l)=1,$ then
$n|l\col(\la) .$ \elem

\prf By monodromy equation $\e(\om^l,\la)\e(\la,\om^l)=
\exp(\frac{2\pi i l \col(\la)}{n})$ and the lemma follows. \qed

\blem\label{shape} If $[v\la]=\sum_{1\leq j\leq {k_1-1}} [\om^{l_1
j} w] $ where $k_1l_1=n, [\om^{jl_1}w]=[\om^{j'l_1}w]$ iff $j=j'\mod
\ k_1,$ and $\sum_{1\leq i\leq n-1}\la_i \leq k-1.$ Then
$\la=(0,...,0,k/k_1,0,...,0,k/k_1,...,0)$ where $(0,...,0,k/k_1)$
(with $l_1-1$ $0$'s)  appears $k_1-1$ times, and the last $l_1-1$
entries are $0$'s, and $\col(\la)=0\mod \  n.$

\elem

\prf Since $[\om^{l_1}\la]=[\la],$ in the  components of $\la,$
$(\la_0,...,\la_{l_1-1})$ appears $k_1$ times.  By assumption $v\la$
is a sum of $k_1$ distinct irreducible subsectors, it follows from
Lemma \ref{fusionrule} that $\la$  has only $k_1$ non-zero
components. Since $\la_0\neq 0,$ and $\col(\la)
=\frac{kl_1(k_1-1)}{2},$ the lemma follows. \qed

\bprop\label{nonprime} If $[a_u]=[x_1y_1], 1<d_{x_1}< d_u $ where
$x_1\prec a_{\la_1}, y_1\prec a_{\la_2}.$ Then either $[x_1]=[a_v],
[y_1]=[b_i]$ or
 $[y_1]=[a_v],
[x_1]=[b_i], 1\leq i\leq n.$ \eprop

\prf By using the action of $\om$ if necessary, we may assume that
the zero-th components of $\la_1,\la_2$ are positive. By Lemma
\ref{klx} we can assume that
\begin{align*}
[a_{\la_1}] &=\sum_{1\leq i\leq k_1} [x_i],
[\om^{l_1}\la_1]=[\la_1], [g^i x_1g^{-i}]=[x_i], 0\leq i\leq k_1-1,
k_1l_1=n \\
[a_{\la_2}]&=\sum_{1\leq i\leq k_2} [y_i], [\om^{l_2}\la_2]=[\la_2],
[g^i x_1g^{-i}]=[x_i], 0\leq i\leq k_2-1,k_2l_2=n
\end{align*}
Since $a_u\prec a_{\la_1\la_2},$ $\col(\la_1) +\col{\la_2}=\col
(u)=1\mod \  n.$ By Lemma \ref{holo} $k_i|\col(\la_i), i=1,2.$ Hence
$(k_1,k_2)=1.$

Since  $x_1y_1, a_\ad$ are irreducible,  we may assume that $\l \bar
x_1x_1, a_\ad\ra=0,$ i.e., $a_vx_1$ is irreducible. Let $w\prec
v\la_1.$ Since $\om^{l_1}[\la_1]=[\la_1],$ $\om^{l_1}w\prec v\la_1.$
Let $t_1|k_1$ be the least positive integer such that
$[\om^{l_1t_1}w]=[w].$ By Lemma \ref{holo} $n|l_1t_1\col(w).$ But
$\col(w)=1+ \col{\la_1}\mod \  n$ with $k_1|\col(\la_1).$ We
conclude that $t_1=k_1$ and
\[[v\la_1]\prec \sum_{0\leq j\leq k_1-1} [ \om^{l_1j} w]\]
Since $a_w\prec a_{v\la_1}=\sum_{1\leq j\leq k_1}[a_v x_j]$ and each
$a_vx_j$ is irreducible, $d_{a_w}=d_w\geq d_vd_{x_1}=
d_vd_{\la_1}/n.$ Hence
\[[v\la_1]= \sum_{0\leq j\leq k_1-1} [ \om^{l_1j} w]\]
By Lemma \ref{shape} we have $\col(\la_1)=0\mod \  n.$ Hence
$\col(\la_2)=1\mod \  n$ and $k_2=1$. If $l_1=1,$ then
$\la_1=(n',...,n'),$ and $d_{\la_2}=d_v.$ By Proposition on Page 10
of \cite{Ganv} $\la_2$ must be in the orbit of $v$ or $\bar v$ under
the action of $\om$. But $\col(\la_2)=1\mod \  n$, so
$[a_{\la_2}]=[a_v]$ and Prop. is proved. In the following we assume
that $l_1\geq 2$ to reach a contradiction. \par

Note that $[a_{\la_1\la_2}]=k_1[a_u]$, hence
$[\la_1\la_2]=\sum_{0\leq i\leq k_1-1} [\om^{l_1i}u].$ By Lemma
\ref{fusionv} $k_1\geq 2.$ We have
\[\l \la_1\la_2,\la_1\la_2\ra= k_1\geq 1+ \l
\la_1\bar\la_1,\ad\ra\l\la_2\bar\la_2,\ad\ra= 1+
(k_1-1)\l\la_2\bar\la_2,\ad\ra\] Hence $\l v\la_2,v\la_2\ra=2.$

On the other hand since $ n=k_1l_1\geq 4$, by Lemma \ref{fusionrule}
we have

$\l \la_1\bar\la_1, (0,1,0,...,0)(0,0,...,1,0)\ra \geq k_1+1,$
$[(0,1,0,...,0)(0,0,...,1,0)]=[v\bar v]+[(0,1,0,...,0,1,0)]$ and we
conclude that
\[\l \la_1\bar\la_1,(0,1,0,...,0,1,0)\ra \geq 1\]
We must have
\[\l(0,1,0,...,0)\la_2, (0,1,0,...,0)\la_2\ra =2\]
Hence by Lemma \ref{fusionrule} $\la_2=(m,0,...,0)$ or
$\la_2=(0,...0,m).$

Note that $[(2,0,...,0)]+[(0,1,0,...,0)]=[v^2].$ If $m>1$ then by
fusion rules
\[
[(2,0,...,0)(0,0,...,2)]= [v\bar v]+ [(2,0,...,2)], \l
(2,0,...,0)\la_2, (2,0,...0)\la_2\ra= 3
\]
We obtain $\l(2,0,...,2), \la_2\bar\la_2\ra =1.$  Similarly we
obtain that  $\l(2,0,...,2), \la_1\bar\la_1\ra \geq 1,$ hence $\l
\la_1\la_2,\la_1\la_2\ra= k_1\geq k_1+1,$ a contradiction. Therefore
$\la_2=v $ or $\bar v.$ But $\col(\la_2)=1\mod \  n$ we have
$\la_2=v.$

From $[\la_1 v]=[\la_1\la_2]=\sum_{0\leq i\leq k_1-1} [\om^{l_1i}u]$
and Lemma \ref{shape} we conclude that $\la_1=(n',n',...,n')$ Hence
$l_1=1$ contradicting our assumption $l_1>1$.\qed

\subsection{The proof of Th. \ref{mn}}\label{finish}
By Lemma \ref{irredad}, Cor. \ref{cor'} and {Prop. \ref{lemh}, the
assumptions of Cor. \ref{smallc} are verified. We can find $\ro_o
c\in H_\ro$ as in Cor. \ref{smallc}. Since
$[\ro_o\bar\ro_o]=\sum_{1\leq i\leq n} [\om^i],$ it follows that
$d_c=1$, and we conclude that $\ro_o c_1\prec \la \ro_o c$ for some
$\la$, and by Prop. \ref{xua} we have
$$1\leq \l \ro_o c_1,
\ro_o a_\la c\ra = \l c_1, \bar\ro_o \ro_o a_\la c\ra=\l c_1, a_\la
\bar\ro_o \ro_o  c\ra
$$
It follows that $ c_1\prec a_\la g^i c$ for some $1\leq i\leq n.$
Since $c_1 (g^ic)^{-1} (M)= c_1(M)$ as a set, replacing $c_1$ by
$c_1(g^ic)^{-1}$ if necessary, we may assume that $[g^ic]=[1],$ and
$c_1\prec a_\la.$ Since $a_u=c_1c_2$ it follows that $c_2\prec
a_\mu$ for some $\mu.$ By Prop. \ref{nonprime} we conclude that
$[c_1]= [a_v], [c_2]=[b_i],$ or $[c_1]=[b_i], [c_2]=[a_v], 1\leq
i\leq n.$ Assume first that $c_1= U a_v U^*, c_2 = U' b_i U'^*$ with
$U, U'$ unitary. Then we have $a_u= ad_{U a_v(U')} a_v b_i =
ad_{U_i} a_v b_i.$ Since $a_v b_i$ is irreducible we have $Ua_v(U')
U_i^*\in \Co,$ and this implies that the intermediate subfactor
$c_1(M)= ad_{U_i} a_v(M)$, i.e, it is one of the subfactors in Prop.
\ref{2n}. The case when $[c_1]=[b_i], [c_2]=[a_v] 1\leq i\leq n.$ is
treated similarly. By Prop. \ref{2n} Th. \ref{mn} is proved.

\section{Related issues}\label{related}
\subsection{Centrality of a class of
intertwinners}\label{centrality} We preserve the general setup of
section \ref{inductionsection}. If $\rho=\mu c, \mu\in \Delta_\A,
d_c=1$ it follows from definition \ref{ala} that $[a_\la]=[\tilde
a_\la] =[ c^{-1}\la c], \forall \la,$ hence $Z_{\la\la_1}=
\delta_{\la,\la_1}.$ Motivated by our proof of Th. \ref{mn} we make
the following:

\bconj\label{conj1} If $Z_{\la\la_1}=\delta_{\la,\la_1},$ then
$\rho=\mu c, \mu\in \Delta_\A, d_c=1.$ \econj

We will prove that Conjecture \ref{conj1} is equivalent to the
centrality of a class of intertwinners. Assume that
$Z_{\la\la_1}=\delta_{\la,\la_1}.$ Then for each irreducible $\la$
there is (up to scalar) a  unique unitary $u_\la\in \Hom
(a_\la,\tilde a_\la).$

Similar to Def. \ref{uwm} we define:
\bdef\label{ula} $u_{\la_1\la_2...\la_m}:= u_{\la_1}
a_{\la_1}(u_{\la_2})... a_{\la_1\la_2...\la_{n-1}}(u_{\la_n})\in
\Hom( a_{\la_1\la_2...\la_m}, \tilde  a_{\la_1\la_2...\la_m})$ \ede

If $\rho=\mu c, \mu\in \Delta_\A, d_c=1,$ then it follows from def.
(\ref{ala}) we can choose $u_\la$ such that $u_\la= c^{-1}(
\tilde\e(\la,\bar\mu) \tilde\e(\bar\mu,\la)).$ Use BFE in Prop.
\ref{bfe} we have
\begin{align*}
u_{\la_1\la_2...\la_m}&= c^{-1}(\tilde\e(\la_1\la_2...\la_m,\bar\mu)
\tilde\e(\bar\mu,\la_1\la_2...\la_m))\in
\Hom(a_{\la_1\la_2...\la_m}, \tilde a_{\la_1\la_2...\la_m}),\\
\Hom(a_{\la_1\la_2...\la_m}, \tilde a_{\la_1\la_2...\la_m})
&=c^{-1}(\Hom (\bar\mu\la_1\la_2...\la_m,
\bar\mu\la_1\la_2...\la_m)).
\end{align*}
By using BFE in Prop. \ref{bfe} again we have proved the following:

\blem\label{muc} If $\rho=\mu c, \mu\in \Delta_\A, d_c=1,$ then
$u_{\la_1\la_2...\la_m} T u_{\la_1\la_2...\la_m}^* =T, \forall T\in
\Hom(a_{\la_1\la_2...\la_m}, a_{\la_1\la_2...\la_m}).$ \elem

Using $u_\la$ we define:
\bdef\label{psi2} For any irreducible $ [b]
\in H_\ro, \la\in \Delta_\A, $
\[
\psi_{b}^{(\la)}:= S_{11} d_bd_\la \phi_\la(
\e(b\bar\ro,\la)b(u_\la)\e(\la,b\bar\ro))
\]
\ede

\blem\label{epsi2}  For any irreducible $ [b] \in H_\ro,
\psi_b^{(\la)}= c_\la \phi_b^{(\la)}, |c_\la c_{\bar\la}| =1$ where
$ c_\la$ are complex numbers independent of $b.$\elem

\prf Since by Lemma \ref{epsi} $\sum_b {\psi_b^{(\la)}}^* [b]$ is an
eigenvector of the action of $\mu$ with eigenvalue
$\frac{S_{\mu\la}}{S_{1\la}},$ and by Prop. \ref{spec} there is up
to scalar a unique such eigenvector, it follows that there is a
complex number $ c_\la$ independent of $b$ such that
$\psi_b^{(\la)}= c_\la \phi_b^{(\la)}, \forall b.$ Similarly  since
$\sum_b {\phi_b^{(\la)}}^* [b]$ is an orthogonal eigenvector of the
action of $\mu$ with eigenvalue
$\frac{S_{\mu\bar\la}}{S_{1\bar\la}},$ we have $\phi_b^{(\bar\la)}=
c_\la' {\phi_b^{(\la)}}^*, |c_\la'|=1, \forall b.$ we have
$\phi_b^{(\bar\la)}= c_{\bar\la} c_\la'{\phi_b^{(\la)}}^*, \forall
b, |c_\la'|=1. $ By Lemma \ref{epsi} $\sum_b
\psi_b^{(\la)}\psi_b^{(\bar\la)}$ has absolute value $1,$ and it
follows that $|c_\la c_{\bar\la}| =1.$ \qed

The following Lemma is proved in the same way as Lemma \ref{indep}:
\blem\label{indep2} If $u_{\la_1\la_2...\la_m}$ is central, then for
fixed $\mu$, if $t_\mu\in \Hom(\mu,\la_1\la_2...\la_m)$ is an
isometry, then
$\bar\ro(t_\mu)^*u_{\la_1\la_2...\la_m}\bar\ro(t_\mu)\in
\Hom(a_\mu,\tilde a_\mu)$ is a unitary independent of the choice of
$t_\mu,$ and is a scalar multiple of $u_\mu$. \elem

\bprop\label{c1=c} Conjecture (\ref{conj1}) is equivalent to the
following statement: If $Z_{\la\la_1}=\delta_{\la,\la_1},$ then
$u_{\la_1\la_2...\la_m}$ is central for all $\la_1,...\la_m, \forall
m.$\eprop

\prf Suppose that Conjecture (\ref{conj1})  is true. Then it follows
from Lemma \ref{muc} that if $Z_{\la\la_1}=\delta_{\la,\la_1},$ then
$u_{\la_1\la_2...\la_m}$ is central for all $\la_1,...\la_m, \forall
m.$ Suppose now that $u_{\la_1\la_2...\la_m}$ is central for all
$\la_1,...\la_m, \forall m.$ As in the proof of Lemma \ref{sumpsi}
by using centrality $u_{\la_1\la_2...\la_m}$  we calculate
\[
\frac{\psi_b^{(\la_1)}}{\psi_b^{(1)}}\frac{\psi_b^{(\la_2)}}{\psi_b^{(1)}}...
\frac{\psi_b^{(\la_m)}}{\psi_b^{(1)}}=\sum_\mu \l
\mu,\la_1...\la_m\ra d_\mu\phi_\mu(\e(b\bar\ro,\mu)b(u_\mu)
\e(\mu,b\bar\ro)) c_\mu
\]
where $|c_\mu|=1.$ Hence using Lemma \ref{epsi} as in the proof of
Lemma \ref{sumpsi} we have
\[
\sum_b
|d_b^2\frac{\psi_b^{(\la_1)}}{\psi_b^{(\la_1)}}\frac{\psi_b^{(\la_2)}}{\psi_b^{(1)}}...
\frac{\psi_b^{(\la_m)}}{\psi_b^{(1)}}|=\l 1,\la_1...\la_m\ra \sum_b
d_b^2 =\sum_\la
\frac{S_{\la_1\la}}{S_{1\la}}\frac{S_{\la_2\la}}{S_{1\la}}...\frac{S_{\la_m\la}}{S_{1\la}}d_\la^2
\]
Now choose $m=2m_1$ and $\la_{i+m_1}=\bar\la_i, 1\leq i\leq m_1,$
sum over $\la_1,...,\la_{m_1}$ and use Lemma \ref{epsi2} we obtain
\[\sum_b \frac{1}{d_b^{m-2}}=\sum_\la \frac{1}{d_\la^{m-2}}
\]
Let $m=2m_1$ go to infinity and notice that $d_b\geq 1$ we conclude
that there must exist a sector $c$ such that $d_c=1$ and $\rho=\mu
c$ for some $\mu\in \Delta_\A$.\qed

%\brem\label{dodd} example shows that centrality is not local issue
%like that of Lemma \ref{325}.\erem

For each irreducible $\la\in \Delta_\A$ we choose $R_{\la\bar\la}$
so that $ R_{\la\bar\la}^* R_{\la\bar\la}= d_\la,
\la(R_{\bar\la\la}^*) R_{\la\bar\la} =1.$ These operatorors are
unique up to scalars.

\blem\label{3251} (1) We can choose $u_\la$ such that
\[\bar\ro(R_{\la\bar\la}^*)u_{\la\bar\la}=\bar\ro(R_{\la\bar\la}^*),
u_{\la\bar\la}\bar\ro(R_{\la\bar\la})
=\bar\ro(R_{\la\bar\la}),\forall \la;
\]
\par
(2) The relative braiding as defined in Lemma \ref{relb} among
$a_\la's$ (resp. $\tilde a_\la$'s )is a braiding and
$\e(a_\la,a_\mu)=\e(\tilde a_\la,\tilde a_\mu)=
\bar\ro(\e(\la,\mu)), \forall \la,\mu\in \Delta_\A.$ \elem

\prf Ad (1): Note that $\bar\ro(R_{\la\bar\la}^*)u_{\la\bar\la}$ is
equal to $\bar\ro(R_{\la\bar\la}^*)$ up to a constant of absolute
value $1$, hence we can choose multiply $u_\la,u_{\bar\la}$ by
suitable constants of absolute value $1$ so that
\[\bar\ro(R_{\la\bar\la}^*)u_{\la\bar\la}=\bar\ro(R_{\la\bar\la}^*)\]
If
\[u_{\la\bar\la}\bar\ro(R_{\la\bar\la})
=c_\la \bar\ro(R_{\la\bar\la}),\forall \la ,\] multiply both sides
on the left by $\bar\ro(R_{\la\bar\la})^*$ we conclude that
$c_\la=1, \forall \la.$\par Ad (2) The relative braidings are
braidings since $[a_\la]=[\tilde a_\la]$ by assumption and Lemma
\ref{relb}. By definition we have
\[ \e(a_\la,a_\mu)=     u_\mu^* \bar\ro(\e(\la,\mu))a_\la(u_\mu)
= u_\mu^* u_\mu \bar\ro(\e(\la,\mu)) =\bar\ro(\e(\la,\mu))
\]
where we have used Lemma \ref{325} in the second $=$ since $u_\mu\in
\Hom(a_\mu,\tilde a_\mu)\subset \Hom(\bar\ro\mu, \bar\ro\mu).$ The
other case is proved similarly.\qed

\bdef\label{cap} An operator is a cap ( resp.cup ) operator if it is
$\mu(R_{\la\bar\la})$ (resp.$\mu(R_{\la\bar\la})^*$) for some
$\mu,\la\in \Delta_\A$. It is a braiding operator it is $
\mu(\e(\la,\nu))$ or $\mu(\tilde\e(\la,\nu))$ for some
$\nu,\mu,\la\in \Delta_\A$.\ede

\bdef\label{bs} Denote by $B_{\la_1\la_2...\la_m}$ the subspace of
$\Hom(\la_1\la_2...\la_m, \la_1\la_2...\la_m)$ which is linearly
spanned by operators in $\Hom(\la_1\la_2...\la_m,
\la_1\la_2...\la_m)$ consisting of products of only caps, cups and
braiding operators. \ede

\bprop\label{ub} For any $T\in \bar\ro(B_{\la_1\la_2...\la_m}),$
$u_{\la_1...\la_m} T = T u_{\la_1...\la_m}.$ \eprop

\prf It is enough to check for an operator $T$ which consists of
products of only caps, cups and braiding operators. Note that the
statement of Prop. is independent of choices of $u_{\la},$ and we
can choose our $u_\la$ so that they verify (1) of Lemma \ref{3251}.
It is useful to think of $T$ as an tangle connecting top $m$ strings
labeled by $a_{\la_1},... a_{\la_m}$ to the bottom $m$ strings
labeled by $a_{\la_1},... a_{\la_m}$ as in Chapter 2 of \cite{Tu},
where in the tangle only cups,caps and braidings are allowed. Then
by Prop. \ref{bfe}, $uTu^*$ will be represented by the same tangle,
except the top and bottom $m$ strings are now labeled by $\tilde
a_{\la_1},... \tilde a_{\la_m}.$  For each  closed string  in
$uTu^*$ labeled by $a_\mu,$ by inserting $u_\mu$ we can change the
label $a_\mu$ to $\tilde a_\mu$ using Prop. \ref{bfe} without
changing the operator since we have a closed  string. Therefore
$uTu^*$ is represented by the same tangle $T$ with all labels
changed from the original labels  $a_\mu$ of $T$ to $\tilde a_\mu.$
Since $T$ consists of products of only caps, cups and braiding
operators, Prop. follows from Lemma \ref{3251}. \qed

\bconj\label{conj2} $B_{\la_1\la_2...\la_m}=
\Hom(\la_1\la_2...\la_m, \la_1\la_2...\la_m), \forall
\la_1,...\la_m, m\geq 1.$\econj

By Prop. \ref{ub} and Prop. \ref{c1=c} we have proved the following:
\bprop\label{c21} Conjecture (\ref{conj2}) implies Conjecture
(\ref{conj1}) . \eprop

By examing the proof of Prop. \ref{c1=c}, we can formulate a weaker
version of Conjecture (\ref{conj2}).

\bdef\label{g} We say that $\la$ is a generator for $\Delta_\A$ if
for any irreducible $\mu\in \Delta_\A,$ there is a positive integer
$m$ such that $\mu\prec \la^m.$ \ede

 \bconj\label{conj3} For some
generator $\la$ of $\Delta_\A$, $B_{\la\la...\la}= \Hom(\la^m,
\la^m), \forall  m\geq 1$ where $m$ is the number of $\la$ that
appears in the definition of $B_{\la\la...\la}.$ \econj

\blem\label{g1} Assume that $\la$ is a generator for $\Delta_\A.$
Then the set $\{[\mu]| |\frac{S_{\la\mu}}{S_{1\mu}}|=d_\la\}$ is a
finite abelian group. \elem

\prf Note that by definition $|\frac{S_{\la\mu}}{S_{1\mu}}|=d_\la$
implies that $\e(\mu,\la)\e(\la,\mu)\in \Co.$ By Prop. \ref{bfe}
this implies that $\e(\mu,\la_1)\e(\la_1,\mu)\in \Co$ if $\la_1\prec
\la^m, m\geq 1$ Since $\la$ a generator, it follows that
 $\e(\mu,\la_1)\e(\la_1,\mu)\in \Co, \forall \la_1\in \Delta_\A. $ Hence
 $|\frac{S_{\mu\la_1}}{S_{1\la_1}}|=d_\mu, \forall \la_1\in \Delta_\A.
 $ By properties of $S$ matrix this implies that $d_\mu =1.$ On the
 other hand if $d_\mu=1$ then $|\frac{S_{\la\mu}}{S_{1\mu}}|=d_\la$
 since $\mu\la$ is irreducible. It follows that the set
 $\{[\mu]| |\frac{S_{\la\mu}}{S_{1\mu}}|=d_\la \}$ is a
finite abelian group.\qed

\bprop\label{c3=c1} Conjecture (\ref{conj3}) implies Conjecture
(\ref{conj1}). \eprop

\prf Assume conjecture (\ref{conj3}) is true. Then by Prop. \ref{ub}
we know that $u_{\la^m}$ is central.

As in the proof of Prop. \ref{c1=c}, replacing $\la_i$ by $\la$ in
the summation we have
\[
|\sum_a (\frac{\psi_a^{(\la)}}{\psi_a^{(1)}})^m d_a^2| = \sum_{\mu}
(\frac{S_{\la\mu}}{S_{1\mu}})^m S_{1\mu}^2
\]
Choose $m$ to be divisible by the order of the finite abelian group
in Lemma \ref{g1} and let $m$ go to infinity, the RHS of the above
equation has leading order (up to multiplication by a positive
number) $d_\la^m.$ It follows that there is a sector $c$ such that
$|\frac{\psi_c^{(\la)}}{\psi_c^{(1)}}|=d_\la.$ For any $\mu\prec
\la^l, l\geq 1.$ Use the centrality of $u_{\la^l}$ we have
\[(\frac{\psi_c^{(\la)}}{\psi_c^{(1)}})^l = \sum_\mu \l \mu,\la^l
\ra \frac{\psi_c^{(\mu)}}{\psi_c^{(1)}} c_\mu\] where $|c_\mu|=1.$
So we have $\sum_{\mu\prec \la^l}
|\frac{\psi_c^{(\mu)}}{\psi_c^{(1)}}| \geq d_\la^l.$ Since
$|\frac{\psi_c^{(\mu)}}{\psi_c^{(1)}}|\leq d_\mu$ and $\sum_{\mu}\l
\mu,\la^l\ra d_\mu =d_\la^l,$ we conclude that
$|\frac{\psi_c^{(\mu)}}{\psi_c^{(1)}}|= d_\mu, \forall \mu\prec
\la^l.$ Since $\la$ is a generator, we conclude that
$|\frac{\psi_c^{(\mu)}}{\psi_c^{(1)}}|= d_\mu, \forall \mu.$ By
Lemma \ref{epsi2} we conclude that
$|\frac{\phi_c^{(\mu)}}{\phi_c^{(1)}}|^2= d_\mu^2.$ Sum over $\mu$
on both sides we conclude that $d_c=1,$ and the Prop. is proved.
\qed

By Prop. \ref{c3=c1} and Lemma \ref{denseb} we have proved the
following:

\bcor\label{hecke} Conjecture (\ref{conj1}) is true for $\Delta_\A$
where $\A$ is the net associated with $SU(n)_k.$ \ecor

\subsection{Maximal subfactors}\label{maxsection}
In this section we give an application of Cor. \ref{hecke}.\par The
following notion is due to V. F. R. Jones: \bdef\label{jones} A
subfactor $N\subset M$ is called maximal if $M_1$ is a von Neumann
algebra such that $N\subset M_1\subset M$ implies $M_1=M$ or
$M_1=N.$\ede

We preserve the setting of section \ref{typea}. We will say that
$\lambda$ is maximal if $\la(M)\subset M$ is a maximal subfactor.

\bprop\label{vnz} If $S_{v\la}\neq 0,$ then $\la$ is maximal. \eprop
\prf Let $M_1$ be an intermediate subfactor between $\la(M)$ and
$M.$  Suppose that $\la=c_1c_2$ and $c_1=c_1'c_1''$ as in Prop.
\ref{loclr}. Since  $S_{v\la}\neq 0,$ apply Lemma \ref{znot0} and
Lemma \ref{loc2} to induction with respect to $c_1'$, we conclude
that $\e(v,c_1'\bar c_1') \e(c_1'\bar c_1',v)\in \Co.$ By Lemma
\ref{scalarv} we conclude that $[c_1'\bar c_1']=[1].$ By Prop.
\ref{loclr} we must have   $Z_{\la1}^{c_1}= \delta_{\la1}.$ Since
$S_{\la v}\neq 0,$ by  Lemma \ref{znot0} and \S2 of \cite{Ganv} we
conclude that $Z_{\mu_1\mu_2}^{c_1}=\delta_{\mu_1\mu_2}.$ By Prop.
\ref{hecke} we conclude that $c_1=\mu c, \mu\in \Delta_\A, d_c=1.$
Replacing $c_1$ by $c_1c^{-1}$ if necessary we may assume that
$c_1=\mu.$ It follows that $c_2=\mu_2$ for some $\mu_2\in
\Delta_\A.$ By Lemma \ref{fusionv} we conclude that $[\mu]=[\la]$ or
$[\mu]=[\omega^i ], 1\leq i\leq n,$ hence $M_1=\la(M)$ or $M_1=M.$
\qed

\bcor\label{vv} If $k+n=p^l$ where $p$ is a prime number, and
$(k,n)\neq (2,2),$ then $\la$ is maximal iff there is no $ 1\leq
i\leq n-1$ such that $[\omega^i\la]=[\la].$\ecor

\prf By Th. 5 of  \cite{GanW} when $k+n=p^l$ where $p$ is a prime
number, $S_{v\a}=0$ iff $[\omega^i\la]=[\la]$ for some $ 1\leq i\leq
n-1.$ Let $i_1|i$ be the smallest positive integer such that
$[\omega^{i_1}\la]=[\la]. $ Then  $[\omega^i\la]=[\la]$ for some $
1\leq i\leq n-1,$ then $[\la\bar\la]\prec \sum_{1\leq j\leq n/i_1}
[\om^{j i_1}]$ and by \cite{ILP} and our assumption that $\la$ is
maximal it follows that $[\la\bar\la]= \sum_{1\leq j\leq n/i_1}
[\om^{j i_1}].$ By Lemma \ref{fusionv} and Lemma \ref{irredad} this
is only possible if $k=n=2.$ The corollary now follows from Prop.
\ref{vnz}.\qed

\blem\label{c12} Assume that $Z_{1\mu}^{c_1}=\delta_{1\mu}, \forall
\mu.$ Then $\langle c_1c_2,c_1c_2\rangle = \langle c_1\bar c_1,\bar
c_2 c_2\rangle.$\elem

\prf By \S2 of \cite{Ganv} we have
$Z_{\mu_1\mu_2}^{c_1}=\delta_{\mu_1\tau(\mu_2)}$ where
$\mu\rightarrow \tau(\mu)$ is an order two automorphism of fusion
algebra. It follows that $[\tilde a_\mu]=[a_\tau(\mu)],$ and by
\cite{BEK2}  irreducible sectors of $\bar c_1 \nu c_1$ are of the
form $a_\mu, \forall \mu.$ Since
$$
\langle c_2\bar c_2, a_\mu\rangle= \langle c_2, a_\mu c_2\rangle=
\langle c_2, c_2\mu\rangle=\langle \bar c_2 c_2,
a_\mu\rangle=\langle a_{\bar c_2 c_2}, a_\mu\rangle,
$$
we conclude that $[c_2\bar c_2]=[a_{\bar c_2 c_2}],$ and
$$
\langle c_1\bar c_1, \bar c_2 c_2\rangle= \langle c_1, \bar c_2 c_2
c_1 \rangle= \langle c_1, c_1 a_{\bar c_2 c_2}\rangle=\langle c_1,
c_1 c_2 \bar c_2\rangle=\langle { c_1 c_2}, c_1 c_2\rangle
$$
\qed

\bcor\label{max} Suppose that $k\neq n-2,n+2,n.$  then $\la$ is
maximal iff there is no $ 1\leq i\leq n-1$ such that
$[\omega^i\la]=[\la].$\ecor

\prf When $k=1$  the Cor. is obvious. By Lemma \ref{irredad} we can
assume that $k\geq 2$ and $d_\ad>1.$
 As in the proof of Cor. \ref{vv}, $\la$ is maximal implies that
there is no $ 1\leq i\leq n-1$ such that $[\omega^i\la]=[\la].$ Now
suppose that there is no $ 1\leq i\leq n-1$ such that
$[\omega^i\la]=[\la].$ If $S_{v\la}\neq 0$, then $\la$ is maximal by
Cor. \ref{vnz}. If $k=2,$ the $S$ matrix elements are equal to that
of $S$ matrix elements for $SU(2)_n$ up to phase factors, and it
follows easily that $S_{v\la}\neq 0$ if there is no $ 1\leq i\leq
n-1$ such that $[\omega^i\la]=[\la].$

Suppose that $k \geq 3, S_{v\la}= 0.$ Since $[v\bar v]=[1]+[\ad]$ we
have $S_{\ad\la}= -S_{1\la}\neq 0.$ Assume that $M_1$ is an
intermediate subfactor between $\la(M)$ and $M$, and $\la=c_1c_2$
with $c_1(M)=M_1$ and $c_1=c_1'c_1''$ as in Prop. \ref{loclr}. Apply
Lemma \ref{znot0} we have $\l a^{c_1'}_{\ad}, \tilde
 a^{c_1'}_{\ad}\ra \geq 1.$ By Lemma \ref{irredad} we must have
 $[ a^{c_1'}_{\ad}]=[ \tilde a^{c_1'}_{\ad}]$ and by Lemma
 \ref{scalarad} $[c_1'\bar c_1']= \sum_{1\leq j\leq n/j_1} [\om^ {j j_1}].$
 By Frobenius reciprocity we have $[\om^{j_1}c_1']=[c_1'.]$
Since $\la=c_1'c_1''c_2,$ $[\om^{j_1}\la]=[\la]$, and by assumption
$j_1=n$ and $[c_1'\bar c_1']=[1]$. By  Prop. \ref{loclr} we must
have   $Z_{\mu1}^{c_1}= \delta_{\mu1}, \forall \mu.$ By \S2 of
\cite{Ganv} we have $Z_{\mu_1\mu_2}^{c_1}=\delta_{\mu_1\tau(\mu_2)}$
where $\tau(\mu)=\omega^{m \col(\mu)} \mu  $ or $\tau(\mu)=\omega^{m
\col(\mu)} \bar\mu,  m\geq 0.$ We claim that in fact
$[\omega^m]=[1]$ and $\tau(\mu)=\mu.$ First we show that
$\tau(\mu)=\omega^{m \col(\mu)} \mu.  $  If instead
$\tau(\mu)=\omega^{m \col(\mu)} \bar\mu,$
 since $k\geq 3,$
$\tau((0,1,0,...,0))\neq (0,1,0,...,0),$ by Lemma \ref{znot0} we
must have $S_{\la(0,1,0,...,0)}=0.$ From the fusion rule
$$
[(0,1,0,...,0)(0,0,...,0,2)]=[(0,1,0,...,0,2)]+[v_0]
$$
we must have $S_{\la(0,1,0,...,0,2)}\neq 0.$ By Lemma \ref{znot0} we
must have $\tau((0,1,0,...,0,2))=(0,1,0,...,0,2)=(2,0,0,...,1,0),$ a
contradiction.  So we conclude that $\tau(\mu)=\omega^{m \col(\mu)}
\mu,\forall \mu$. It follows that $[\tilde a_\mu]=[a_{\omega^{m
\col(\mu)}}a_\mu],$ and in particular $[\tilde
a_v]=[a_{\omega^{m}}a_v]$. So we have
$$
[\omega^m v c_1]=[c_1 \tilde a_v]=[c_1 a_v]=[vc_1],
$$
and similarly $[c_2 \omega^{-m}\bar v ]=[c_2\bar v]$. If
$[\omega^m]\neq [1],$ by our assumption on $\la$ we have
$\omega^m\not\prec c_1\bar c_1, \omega^m\not\prec \bar c_2 c_2 .$ On
the other hand we have
$$
\langle \bar v \omega^m v,c_1\bar c_1\rangle \geq 1, \langle \bar v
\omega^m v,\bar c_2 c_2\rangle \geq 1
$$
It follows that $\omega^m v_0\prec c_1\bar c_1, \omega^m v_0\prec
\bar c_2 c_2,$ and $\langle c_1\bar c_1, \bar c_2 c_2\rangle \geq
2.$ By Lemma \ref{c12} we conclude that $\la=c_1c_2$ is not
irreducible, contradicting our assumption. Hence $[\omega^m]=[1]$
and $Z_{\mu_1\mu_2}=\delta_{\mu_1\mu_2}.$  The rest of the proof now
follows in exactly the same way as in the proof of Prop.
\ref{vnz}.\qed
\begin{example}
When $n=2$ we have Jones subfactors and their reduced subfactors. In
the case $k=n=2$ there are three irreducible subfactors and they are
maximal. Let $n=2, k\neq 2.$ Then $\la$ can be labeled by an integer
$1\leq i\leq k.$ Cor. \ref{max} implies that $i$ is maximal iff
$i\neq k/2$ (When $k=4$ this can be easily checked directly). This
can also be proved directly using the same argument at the end of
section \ref{orb}.
\end{example}

\end{document}